\documentclass[11pt]{article}
\usepackage{amssymb, amsmath, amscd, amsthm}
\usepackage[all]{xy}

\usepackage[OT2,OT1]{fontenc}

       {\arraycolsep 0.14em\begin{eqnarray}}{\end{eqnarray}}
\newenvironment{Eqnarray*}%
       {\arraycolsep 0.14em\begin{eqnarray*}}{\end{eqnarray*}}
       {\arraycolsep 0.14em\begin{array}}{\end{array}}

\theoremstyle{plain}
\newtheorem{Thm}{Theorem}[section]
\newtheorem{Lem}[Thm]{Lemma}
\newtheorem{Prop}[Thm]{Proposition}
\newtheorem{Cor}[Thm]{Corollary}

\newtheorem{Conj}[Thm]{Conjecture}

\theoremstyle{definition}
\newtheorem{Def}[Thm]{Definition}

\theoremstyle{remark}
\newtheorem{Rem}[Thm]{Remark}

\numberwithin{equation}{section}

\def\action{\, {\scriptscriptstyle \stackrel{\circ}{{}}} \, }

  \def\C{\mathbb C}  \def\R{\mathbb R}
\def\Z{\mathbb Z}    \def\bP{\mathbb P}
 \def\bD{\mathbb D}

\def\a{\alpha}  \def\b{\beta}     \def\g{\gamma}   \def\d{\delta}
\def\c{\theta}  \def\l{\lambda}   \def\p{\phi}     \def\s{\sigma}
\def\t{\tau}    \def\vp{\varphi}  \def\r{\rho}     \def\y{\eta}
\def\z{\zeta}   \def\w{\omega}    \def\e{\varepsilon}
\def\x{\xi}       

\def\G{\Gamma}  \def\L{\Lambda}      

      \def\cF{\mathcal F}
      
     \def\cN{\mathcal N}
\def\cO{\mathcal O}      
      
 \def\cX{\mathcal X}  \def\cZ{\mathcal Z}  

\DeclareMathAlphabet{\mathzap}{OT1}{pzc}{m}{it}
 \def\zD{\mathzap D}   \def\zE{\mathzap E}

\def\cyr{%
\renewcommand\rmdefault{wncyr}%
\renewcommand\sfdefault{wncyss}%
\renewcommand\encodingdefault{OT2}%
\normalfont
\selectfont}
\DeclareTextFontCommand{\textcyr}{\cyr}
\newcommand{\Dye}{\mbox{\cyr D}}
\newcommand{\Lye}{\mbox{\cyr L}}
\newcommand{\pye}{\mbox{\cyr p}}

\newcommand{\xye}{\mbox{\cyr x}}

\def\Gf{\mathfrak f}       
  \def\Gm{\mathfrak m}     
\def\Gp{\mathfrak p}      

\def\GB{\mathfrak B}
\def\GC{\mathfrak C}  \def\GD{\mathfrak D}  \def\GF{\mathfrak F}
\def\GS{\mathfrak S}   \def\GU{\mathfrak U}  \def\GV{\mathfrak V}
\def\GW{\mathfrak W}

     \def\PSL{\operatorname{PSL}}      
\def\PSU{\operatorname{PSU}}

\def\Aut{\operatorname{Aut}}

  \def\HH{\operatorname{H}}

\def\Span{\operatorname{Span}}

\def\Imag{\operatorname{Im}}    \def\Real{\operatorname{Re}}      

\def\RP{\R\bP}
\def\CP{\C\bP} 

\def\del{\partial}

\def\({ \left( }     \def\){ \right) }
\def\<{ \left\langle } \def\>{ \right\rangle }
\def\hasira{\rule{0mm}{2eX}}

\newcommand{\pd}[1]{\frac{\del}{\del #1}}
\newcounter{mynum}

\title{A construction of Einstein-Weyl spaces \\ 
 via LeBrun-Mason type twistor correspondence}

\author{Fuminori Nakata\thanks{This work 
 is partially supported by Grant-in-Aid for Scientific Resaerch 
 of the Japan Society for the Promotion of Science. }
}

\date{June 17, 2008}

\begin{document}
\maketitle

\abstract{We construct infinitely many Einstein-Weyl structures
on $S^2\times\R$ of signature $(-++)$ which is sufficiently close to the model case of
constant curvature, and whose space-like geodesics are all closed.
Such structures are obtained from small perturbations of the diagonal of
$\CP^1\times\overline{\CP^1}$ using the method of LeBrun-Mason type
twistor theory. The geometry of constructed Einstein-Weyl space 
is well understood from the configuration of holomorphic discs. 
We also review Einstein-Weyl structures and their properties 
in the former half of this article.}

\vspace{8mm}
{\noindent
 {{\it Mathematics Subject Classifications} (2000) : 53C28, 32G10, 53C50, 
  53A30, 53C25. }

\section{Introduction}

Twistor type correspondences for the following structures are known
(See \cite{bib:Hitchin82}): 
\begin{list}{(T\,\arabic{mynum})}%
  {\usecounter{mynum} \itemsep 0in  \leftmargin .4in} 
 \item projective structures on complex 2-manifolds, 
 \item self-dual conformal structures on complex 4-manifolds, 
 \item Einstein-Weyl structures on complex 3-manifolds.
\end{list}
(T2) is the original twistor theory introduced by R.~Penrose \cite{bib:Penrose76}. 
(T3) is called Hitchin correspondence or mini-twistor correspondence. 

There are many progresses on these twistor theory; 
more detail or concrete investigation \cite{bib:Pedersen86,bib:PT93}, 
real objects and reduction theory 
\cite{bib:Calderbank,bib:DW,bib:DWII,bib:JT,bib:Tod92}, 
relation with the theory of integrable systems 
\cite{bib:Dunajski03,bib:DMT00}, and so on. 
The geometric structures treated in these literature are either 
complex or real slices of complex objects, hence they are all analytic. 

On the other hand, real indefinite case, for example, admits 
non-analytic solutions. 
Recently, C.~LeBrun and L.~J.~Mason developed another type of 
twistor theory by which we can also treat such non-analytic solutions 
\cite{bib:LM02,bib:LM05} (see also \cite{bib:Nakata,bib:NakataII}). 
The structures investigated by LeBrun and Mason are 
\begin{list}{(LM\,\arabic{mynum})}%
  {\usecounter{mynum} \itemsep 0in  \leftmargin .4in}
 \item Zoll projective structures on $S^2$ or $S^2/\Z_2$, and 
 \item self-dual conformal structures of signature $(++--)$ on $S^2\times S^2$
   or $(S^2\times S^2)/\Z_2$. 
\end{list}
Here, a projective structure is called Zoll if and only if 
all the maximal geodesics are closed. 
Notice that (LM1), (LM2) are the real objects corresponding to (T1), (T2) respectively. 

There are several remarkable points for LeBrun-Mason theory. 
First, the twistor space is given as a pair $(Z,N)$ of a complex manifold $Z$ and 
a totally real submanifold $N$ in $Z$. 
The ``twistor lines", or in other words, 
the ``nonlinear gravitons" are given by holomorphic disks on $Z$ whose boundaries 
lie on $N$ while the twistor lines are embedded $\CP^1$ in 
Penrose's or Hitchin's case. 
Second, the structures (LM1) and (LM2) are obtained from a small perturbations 
of $N$ in $Z$. So we can treat only the objects which are sufficiently close to 
the model case so far. 
For the last, the corresponding geometry satisfies a global condition, 
for example, Zoll condition in (LM1) case. 

Based on these backgrounds, in this article, 
we investigate in another possibility,
the LeBrun-Mason type correspondence for Einstein-Weyl structures.  
Let us review the definitions and then we state the conjecture and the main theorem. 
Let $X$ be a real (or complex) manifold. 
\begin{Def}
 Let $[g]$ be the conformal class of a definite or an indefinite metric $g$ 
 (or holomorphic bilinear metric for the complex case) on $X$, 
 and $\nabla$ be a (holomorphic) connection on $TX$. 
 The pair $([g],\nabla)$ is called {\it Weyl structure} on $X$ if 
 there exist a (holomorphic) 1-form $a$ on $X$ such that  
 \begin{equation} \label{eq:compatibility_of_Weyl_str}
  \nabla g = a\otimes g. 
 \end{equation}  
\end{Def}
\begin{Def}
 A Weyl structure $([g],\nabla)$ is called {\it Einstein-Weyl} if the symmetrized 
 Ricci tensor $R_{(ij)}=\frac{1}{2}(R_{ij}+R_{ji})$ is proportional to the metric 
 tensor $g_{ij}$, i.e. if we can write 
 \begin{equation}
  R_{(ij)}=\L\, g_{ij}
 \end{equation}
 using a function $\L$ which depends on the choice of $g\in[g]$. 
\end{Def}

Let $[g]$ be an indefinite conformal structure on a real manifold $X$. 
A tangent vector $v$ on $X$ is called {\it time-like} if $g(v,v)<0$, 
{\it space-like} if $g(v,v)>0$ and {\it light-like} or null if $g(v,v)=0$. 
We introduce the following global condition. 
\begin{Def}
 An indefinite Weyl structure $([g],\nabla)$  is called {\it space-like Zoll} 
 if and only if every maximal space-like geodesic is closed.  
\end{Def}

Now we state the conjecture for the LeBrun-Mason type correspondence 
for Einstein-Weyl structures. 
\begin{Conj}
 There is a natural one-to-one correspondence between 
 \begin{itemize}
  \item equivalence classes of 
    space-like Zoll Einstein-Weyl structures on $S^2\times \R$;  and
  \item  equivalence classes of totally real embeddings 
    $\iota : \CP^1\hookrightarrow \CP^1\times \CP^1$, 
 \end{itemize} 
 at least in a neighborhood of the standard objects. 
\end{Conj}
Here the standard embedding $\CP^1\hookrightarrow\CP^1\times\CP^1$ 
is given by
 $\z\mapsto(\z,\bar{\z}^{-1})$ using the inhomogeneous coordinate 
 $\z$ of $\CP^1$. 
The standard Einstein-Weyl structure is 
explained in Section \ref{Section:The_standard_case}. 
Before we state the main theorem, we define the following notion. 
\begin{Def} \label{Def:adapted_vector}
 Let $Z$ be a complex manifold and $\GD\subset Z$ be a holomorphic 
 disk whose boundary is embedded in $Z$. Let $v\in T_pZ$ be a 
 non zero tangent vector at $p\in\del\GD$. 
 Then $v$ is called to be adapted to $\GD$ and denoted by 
 $v\parallel\GD$ iff $v\in T_p\del\GD$ and $v$ has the 
 same orientation as the orientation of $\del\GD$ 
 which is induced from the complex orientation of $\GD$. 
\end{Def}

The main theorem is the following. This 
is the half of the correspondence in the above conjecture; 
from the embedding $\iota$ to the Einstein-Weyl space. 
We also claim that the geometry of the constructed Einstein-Weyl space 
is characterized by the holomorphic disks in the following way. 

\begin{Thm} \label{thm:main_theorem}
 Let $N$ be the image of any embedding of $\CP^1$ into $Z=\CP^1\times\CP^1$ which 
 is $C^{2k+5}$ close to the standard one. 
 Then there is a unique family of holomorphic disks 
 $\{\GD_x\}_{x\in S^2\times\R}$ 
 such that each boundary $\del\GD_x$ lies on $N$, and that 
 the parameter space $M=S^2\times\R$
 has a unique $C^k$ indefinite Einstein-Weyl structure $([g],\nabla)$ 
 satisfying the following properties.  
 \begin{enumerate}
  \item For each $p\in N$, 
    $\GS_p= \left\{x\in M \, | \, p\in \del \GD_x \right\}$ 
	is maximal connected null surface on $M$ 
	and every null surface can be written in this form. 
  \item For each $p\in Z\setminus N$, 
    $\GC_p=\left\{x\in M \, | \, p\in \GD_x \right\}$ 
    is maximal connected time-like geodesic and every 
	time-like geodesic on $M$ can be written in this form. 
  \item For each $p\in N$ and non zero $v\in T_pN$, 
    $\GC_{p,v}= \left\{x\in M \, | \, p\in \del \GD_x, v\parallel \GD_x \right\}$ is 
	maximal connected null geodesic on $M$ and every 
	null geodesic on $M$ can be written in this form. 
  \item For each distinguished $p,q\in N$, 
    $\GC_{p,q}= \left\{x\in M \, | \, p, q \in \del \GD_x \right\}$ is 
	connected closed space-like geodesic on $M$ and every 
	space-like geodesic on $M$ can be written in this form. 
 \end{enumerate}
 In particular, this Einstein-Weyl structure is space-like Zoll. 
\end{Thm}

The organization of this paper is the following. 
We first review the projective structure in Section \ref{Section:Projective_structures}. 
Next, we study about Einstein-Weyl spaces of 
complex, definite or indefinite cases separately
Section \ref{Section:Einstein-Weyl_structures}. 
We prove that, for each case, the Einstein-Weyl condition can be translated to an 
integrability condition for certain distributions. 
Applying this method, we review the proof of Hitchin correspondence 
in Section \ref{Section:Hitchin_correspondence}. 
In Section \ref{Section:The_standard_case}, 
the model case of the LeBrun-Mason type correspondence is explained. 
The standard Einstein-Weyl space is obtained 
as a double cover of a real slice of the Hitchin's example. 
We also study about detail properties for this model case. 

From Section \ref{Section:Holomorphic_disks}, 
we treat the perturbation of the model case. 
In Section \ref{Section:Holomorphic_disks}, we prove that, 
for a small perturbation of the real submanifold $N$, there is a unique 
family of holomorphic disks whose boundaries lie on $N$. 
This family preserves similar properties as the model case, 
especially concerning the double fibration, 
which is studied in Section \ref{Section:The_double_fibration}. 
Finally in Section \ref{Section:Construction_of_E-W_spaces}, 
we prove that there is a unique Einstein-Weyl structure on the parameter space 
of the constructed family of holomorphic disks. 
We also prove that the geometry of the Einstein-Weyl space is 
characterized by the holomorphic disks as in Theorem \ref{thm:main_theorem}.

\section{Projective structures}
\label{Section:Projective_structures}

In this section, we review the projective structures. 
Let $X$ be a real smooth 
$n$-manifold and $x^i\, (i=1, \cdots, n)$  be a local coordinate on $X$. 
The following arguments also works well in the complex case 
by considering $x^i$ as a complex 
coordinate, and using holomorphic functions instead of smooth functions. 

\begin{Def}
 Two connections $\nabla$ and $\nabla'$ on the tangent bundle $TX$ are called 
 {\it projectively equivalent} if their geodesics coincide without considering 
 parameterizations. A projectively equivalent class $[\nabla]$ is called  
 a {\it projective structure} on $X$. 
\end{Def}

Let $\nabla$ and $\nabla'$ be connections on $TX$, and let  
$\G^i_{jk}$ and ${\G'}^i_{jk}$ be their Christoffel symbols respectively, i.e. 
$ \nabla_{\del_k}\del_j=\sum \G^i_{jk} \del_i $
and so on, where we denote $\del_i=\pd{x^i}$. 
Notice that $\nabla$ is torsion-free if and only if $\G^i_{jk}=\G^i_{kj}$. 

\begin{Prop}
 Suppose that both $\nabla$ and $\nabla'$ are torsion-free, then they  
 are projectively equivalent if and only if there exist 
 functions $f_i\, (i=1,\cdots,n)$ on $X$ and the following condition holds: 
 \begin{equation} \label{eq:proj_eq_in_Gamma}
   \G^i_{jk}={\G'}^i_{jk}+\frac{1}{2}(\d^i_j f_k+\d^i_k f_j). 
 \end{equation}
\end{Prop}
\begin{proof}
 A curve $\g: (-\e,\e)\to X: t \mapsto (x^i(t)) $ is a geodesic for $\nabla$
 if and only if 
 $$ \frac{d^2x^i}{dt^2}+\G^i_{jk}\frac{dx^j}{dt}\frac{dx^k}{dt}=0, $$
 which is equivalent to the following equations  
 \begin{equation} \label{eq:geodesic_sys}
  \frac{dx^i}{dt}=y^i, \quad \frac{dy^i}{dt}=-\G^i_{jk}y^j y^k. 
 \end{equation}
 Notice that the natural lifting 
 $\tilde{\g}: t \mapsto (\g(t),\g'(t))\in TX$ of $\g$ is an
 integral curve of the vector field 
 \begin{equation} \label{eq:vf_of_geodesic_spray}
  v= y^i\pd{x^i} - \G^i_{jk} y^j y^k \pd{y^i} 
 \end{equation}
 on $TX$, where $(y^i)$ is the fiber coordinate on $TX$ with respect to 
 the frame $\{\pd{x^i}\}$. 
 Let $\pi:TX\setminus 0_X \to \bP(TX)$ be the projectivization 
 where $0_X$ is the zero section. 
 Then $\pi_*(v)$ defines a line distribution on $\bP(TX)$ 
 which is the {\it geodesic spray}, 
 the distribution defined from the natural lifts of the geodesics. 

 Let $v'$ be the vector field on $TX$ obtained from $\nabla'$ 
 similar as (\ref{eq:vf_of_geodesic_spray}). 
 Notice that $v$ and $v'$ induces the same distribution on $\bP(TX)$ 
 iff $v-v'$ is proportional to $\sum y^i\pd{y^i}$, i.e. 
 \begin{equation} \label{eq:difference_of_G_and_G'}
  \G^i_{jk} y^j y^k \pd{y^i}= {\G'}^i_{jk} y^j y^k \pd{y^i} + f y^i\pd{y^i},  
 \end{equation}
 for some function $f$ on $TX$. 
 Comparing each side, $f$ must be degree one polynomial, so we can write $f=f_iy^i$. 
 Then (\ref{eq:difference_of_G_and_G'}) is equivalent to (\ref{eq:proj_eq_in_Gamma}) 
 since we have $\G^i_{jk}=\G^i_{kj}$ and $\G'^i_{jk}=\G'_{kj}$ from the 
 torsion free condition. 
 Since $\nabla$ and $\nabla'$ are projectively equivalent if and only if 
 the geodesic sprays coincide, we obtain the statement. 
\end{proof}

\begin{Rem} \label{rem:constructing_connection}
 Let $G^i \, (i=1,\cdots,n)$ 
 be functions on $TX$ each of which is a degree-two polynomial for $y^i$. 
 Then, the vector field 
 $$ v= y^i\pd{x^i}-G^i\pd{y^i} $$ 
 on $TX$ defines a torsion-free connection by defining $\G^i_{jk}$ by 
 $G^i=\G^i_{jk}y^j y^k$. 
\end{Rem}

In the complex case, we can prove the following. 
\begin{Prop} \label{prop:family_of_curves=>proj_str}
 Let $X$ be a complex $n$-manifold, and $\cF$ be a holomorphic family of 
 holomorphic curves on $X$. Suppose that, for each non-zero 
 tangent vector $v\in TX$, 
 there is a unique member of $\cF$ which tangents to $v$. 
 Then there is a unique projective structure $[\nabla]$ on $X$ so that  
 $\cF$ coincides to the family of geodesics. 
\end{Prop}
\begin{proof}
 Let $TX\setminus 0_X \overset{\pi}\longrightarrow \bP(TX) 
 \overset{p}\longrightarrow X$ be the projections, 
 $(x^i)$ be a local coordinate on $X$, 
 and $(y^i)$ be the fiber coordinate with respect to the frame $(\pd{x^i})$. 
 A holomorphic curve $c$ on $X$ lifts canonically to the curve $\tilde{c}$ on $\bP(TX)$, 
 and the velocity vector field of $\tilde{c}$ extends to a vector field 
 on $\pi^{-1}(\tilde{c})\subset TX$ of the form 
 \begin{equation} \label{eq:vector_field_induced_from_cF}
  v=y^i\pd{x^i} + G^i(y)\pd{y^i} 
 \end{equation}
 where $G^i(y)$ is a function on $\pi^{-1}(\tilde{c})$ with 
 homogeneity $2$ for $y$, i.e. $G^i(ay)=a^2 G^i(y)$ for every 
 $a\in\C^\times$. 
 
 Since the statement is local, we can assume $\bP(TX)=X\times\CP^{n-1}$. 
 Let $\CP^{n-1}=\cup W_\a $ be an affine open cover. 
 Applying the above method to the curves of $\cF$, 
 we obtain a holomorphic vector field on 
 each $X\times W_\a$ of the form 
  \begin{equation}
  v_\a=y^i\pd{x^i} + G^i_\a(y)\pd{y^i}, 
 \end{equation}
 where $G^i_\a(y)$ is a holomorphic function on $X\times W_\a$ with
 homogeneity $2$ for $y$. Since $v_\a$ and $v_\b$ induce the same 
 geodesic spray, we can write 
 $v_\a-v_\b=f_{\a\b}(y)y^i\pd{y^i}$ on $X\times W_\a \cap X\times W_\b$
 using a holomorphic function $f_{\a\b}(y)$ of homogeneity $1$ for $y$. 

 Since $\HH^1(\bP^{n-1},\cO(1))=0$, we can take $\{v_\a\}$ 
 so that $f_{\a\b}=0$. 
 Hence we obtain a vector field on whole $\bP(TX)$ of the form 
 (\ref{eq:vector_field_induced_from_cF}). 
 Then $G^i$ must be a degree-two polynomial, so we obtain a torsion-free connection 
 $\nabla$ by putting $G^i(y)=\G^i_{jk}y^jy^k$ 
 (Remark \ref{rem:constructing_connection}). 
 Here $\nabla$ is determined up to projective equivalence
 since the ambiguity of taking $v$ remains. 
\end{proof}

\section{Einstein-Weyl structures}
\label{Section:Einstein-Weyl_structures}

In this section, we study about the basic properties of 
3-dimensional Einstein-Weyl structures. 
We will prove that the Einstein-Weyl condition 
is equivalent to the integrability condition of certain distributions. 
We treat three cases separately, i.e. complex, definite, 
and indefinite cases. \\

\noindent
{\bf complex case :}
Let $X$ be a complex 3-manifold and  
$([g],\nabla)$ be a Weyl structure on $X$. 
Though we argue for fixed $g\in[g]$, the statements 
do not depend on the choice of $g$. 
We denote 
$$ \begin{aligned}
 T_\C X &= \, TX\otimes\C \, =T^{1,0}X\oplus T^{0,1}X, \\ 
 T^*_\C X&=T^*X\otimes\C=T^{*\,1,0}X\oplus T^{*\, 0,1}X. 
 \end{aligned} $$
Notice that $g$ induces complex bilinear metrics on 
$T^{1,0}X, T^{0,1}X, T^{*\, 1,0}X$, and $T^{*\, 0,1}X$ which we also denote $g$. 

\begin{Def}
 For each $x\in X$, a complex two-dimensional subspace 
 $V\subset T_x^{1,0} X$ 
 is called {\it null plane} if the restriction of $g$ on $V$ degenerates. 
\end{Def}
The following property is easily checked. 
\begin{Lem}
 If $v\in T_x^{1,0}X$ is null, then $v^\perp$ is a null plane. 
 Conversely, every null plane is written as $v^\perp$ for some null vector $v$. 	 
\end{Lem}

Notice that $v^\perp =\ker v^*$ for every $v\in T_x^{1,0} X$ where 
$v^*=g(v,\cdot) \in T_x^{*\, 1,0}X$, 
and that $v$ is null if and only if $v^*$ is null.  
Let $N(T^{*\, 1,0} X)$ be the null cotangent vectors, 
and $\cZ=\bP(N(T^{*\, 1,0}X))$ be its complex projectivization. 
Notice that each point $u\in\cZ$ corresponds to the null plane 
$V_u=\ker\l$ where 
$\l\in N(T^{*\, 1,0} X)$ is the cotangent vector satisfying $u=[\l]$. 
We can define a complex $2$-plane distribution 
$\zD\subset T^{1,0}\cZ$ 
so that $\zD_u\subset T^{1,0}_u\cZ$ is the horizontal lift of the null plane 
$V_u$ with respect to $\nabla$. 
Notice that the horizontal lift is well-defined 
since $N(T^{*\, 1,0} X)$ is parallel to $\nabla$ because of the 
compatibility condition (\ref{eq:compatibility_of_Weyl_str}).

\begin{Prop} \label{prop:integrability_of_complex_case}
 Let $X$ be a complex 3-manifold. 
 A Weyl structure $([g],\nabla)$ with torsion-free $\nabla$ on $X$ 
 is Einstein-Weyl if and only if 
 the induced distribution $\zD$ on $\cZ$ is integrable, i.e. involutive. 
\end{Prop}
\begin{proof}
 Let $\{e_1,e_2,e_3\}$ be an orthonormal complex local frame on $T^{1,0} X$ 
 with respect to $g\in[g]$, and $\{e^1,e^2,e^3\}$ 
 be the dual frame on $T^{*\, 1,0}X$. 
 Let $\w=(\w^i_j)$ be the connection form of $\nabla$ with respect to $\{e_i\}$, 
 and let $K^i_j=K^i_{jkl}e^k\wedge e^l$ be its curvature form. 
 Then from the compatibility condition (\ref{eq:compatibility_of_Weyl_str}), 
 we obtain the following symmetry for $K$ : 
 \begin{equation} \label{eq:symmetries_of_K}
  \begin{array}{c} K^i_{jkl}=A^i_{jkl}+\d^i_jB_{kl}, \\[1.5eX]
    A^i_{jkl}=-A^i_{jlk}=-A^j_{ikl}, \quad B_{kl}=-B_{lk}. 
  \end{array}
 \end{equation} 
 Since the frame is orthonormal, Einstein-Weyl equation is 
 \begin{equation*} 
   R_{(12)}=R_{(23)}=R_{(31)}=0, \quad R_{(11)}=R_{(22)}=R_{(33)}, 
 \end{equation*}
 and this is equivalent to 
 \begin{equation} \label{eq:EW_eq_in_A}
   A^1_{213}+A^1_{312}=A^2_{321}+A^2_{123}=A^3_{132}+A^3_{231}=0, 
  \quad A^1_{212}=A^2_{323}=A^3_{131}.
 \end{equation}
 
 Now let $\cN=N(T^{*\, 1,0}X)\setminus 0_X$, and $\pi:\cN\rightarrow\cZ$ 
 be the projection where $0_X$ is the zero section. 
 Then $\zD$ is integrable if and only if the pull-back 
 $\pi^*\zD$ is integrable. 
 Here $\pi^*\zD\subset T^{1,0}\cN$ 
 is the complex $3$-plane distribution defined by 
 $\pi^*\zD=\{v\in T\cN \, | \, \pi_*(v)\in \zD\}$. 
 On the other hand, there is a $2$-plane distribution 
 $\tilde{\zD}\subset T^{1,0}\cN$ 
 which is defined in the similar way to $\zD$, i.e. 
 $\zD_u$ is the horizontal lift of the null plane $V_u$. 
 These distributions are related by 
 $\pi^*\zD= \tilde{\zD} \oplus \<\Upsilon\> $
 where 
 \begin{equation} \label{eq:Upsilon}
  \Upsilon=\sum \l_i\pd{\l_i}
 \end{equation}
 is the Euler differential. 
 Now we define several 1-forms on $\cN$ by 
 \begin{equation} \label{eq:c&c_i&t_ij}
  \c=\sum\l_i e^i, \quad \c_i=d\l_i-\sum \l_j\w^j_i, \quad 
  \t_{ij}=\l_i \c_j - \l_j \c_i. 
 \end{equation}
 Then $\tilde{\zD}=\{v \in T\cN \, |\, \c(v)=\c_i(v)=0\, (\forall i)\, \}$ 
 and  $\pi^*\zD=\{v \in T\cN \, |\, \c(v)=\t_{ij}(v)=0 \, (\forall i,j)\, \}$. 
 Hence $\zD$ is integrable if and only if the 1-forms $\{\c, \t_{ij}\}$ on $\cN$ 
 are involutive. Notice that 
 $\t_{23}/\l_1=\t_{31}/\l_2=\t_{12}/\l_3$, hence 
 $\t_{ij}$ are proportional to each other. 

 Let us prove that $\zD$ is integrable if and only if (\ref{eq:EW_eq_in_A}) holds. 
 First, we claim that $d\c\equiv 0 \mod \<\c, \t_{ij}\>$ always holds. 
 Indeed, since $\c_1/\l_1\equiv\c_2/\l_2\equiv\c_3/\l_3$, we have 
 $$ \sum \c_i\wedge e^i \equiv \frac{\c_1}{\l_1}\wedge\c \equiv 0 \ 
 \mod  \<\c, \t_{ij}\>. $$ 
 On the other hand, we have the torsion-free condition: 
 $de^i+\sum\w^i_j\wedge e^j=0$.  
 Then 
 $$ d\c=\sum d\l_i\wedge e^i + \sum \l_i de^i
  = \sum \c_i\wedge e^i + \sum \l_i(de^i+\w^i_j\wedge e^j) 
  \equiv 0 \ \mod \<\c,\t_{ij}\>. $$  
 Next, a direct calculation shows that 
 \begin{equation}
  d\t_{12}\equiv -\sum\l_1\l_j K^j_2 + \sum\l_2\l_j K^j_1 \quad \mod \t_{12}, 
 \end{equation} 
 and we can check that $d\t_{12}\equiv 0$ holds if and only if 
 \begin{equation*}
  \begin{aligned}
  0 = \l_3 \left[  \hasira \right. &
       -A^2_{323}\l_1^2 - A^3_{131} \l_2^2 - A^1_{212}\l_3^2 \\
   & \left. \hasira + (A^3_{132}+A^3_{231})\l_1\l_2 
     + (A^2_{321}+A^2_{123})\l_3\l_1 + (A^1_{213}+A^1_{312})\l_2\l_3 \right]. 
  \end{aligned}
 \end{equation*}
 for every $(\l_i)$ satisfying $\sum\l_i^2=0$. Hence $\zD$ is integrable 
 if and only if the Einstein-Weyl equation (\ref{eq:EW_eq_in_A}) holds. 
\end{proof}

$\zD$ can be explicitly described in the following way. 
As in the above proof, let us 
take a local orthonormal frame $\{e_1, e_2, e_3\}$ on an open set $U\subset X$. 
From the compatibility condition (\ref{eq:compatibility_of_Weyl_str}), 
the connection form $\w$ of $\nabla$ is written 
\begin{equation} \label{eq:conn_form_cpx}
 \w=\begin{pmatrix} 
      \p & \y^1_2 & \y^1_3 \\[1mm]
      \y^2_1 & \p & \y^2_3 \\[1mm]
	  \y^3_1 & \y^3_2 & \p \end{pmatrix}, 
 \quad \text{with} \ 
 \y^j_i=-\y^i_j. 
\end{equation}
We can write 
$$ \begin{aligned}
  N(T^{*\, 1,0}X)|_U &=\left\{ \left. 
    \sum \l_i e^i \, \right| \, \sum\l_i^2=0 \right\}, \\ 
  \cZ|_U &=\left\{ [\l_1:\l_2:\l_3] \, \left| \, \sum\l_i^2=0 \right. \right\}. 
 \end{aligned} $$
Then we obtain 
\begin{equation}
\t_{23}= \l_2d\l_3- \l_3d\l_2 
  + \l_1 \(\l_1\y^2_3 + \l_2\y^3_1 + \l_3\y^1_2  \). 
\end{equation}

Let $U \times \CP^1 \overset\sim\to \cZ|_U$ be a trivialization given by 
\begin{equation}
 (x,\z) \longmapsto [i(1+\z^2):1-\z^2:2\z]
\end{equation}
where $\z\in\C\cup\{\infty\}$ is a inhomogeneous coordinate. 
The horizontal lift $\tilde{v}$ of $v\in T_xU$ at $(x,\z)\in\cZ|_U$ is 
\begin{equation} \label{eq:horizontal_lift_cpx}
 \tilde{v}= v + \left\{ \frac{\y^2_3+i\y^1_3}{2} - i\z\y^1_2 
    + \z^2 \frac{\y^2_3-i\y^1_3}{2} \right\}(v) \pd{\z}. 
\end{equation}
For $(x,\z)\in\cZ|_U$, the corresponding null plane on $T^{1,0}_xX$ is spanned by 
\begin{equation} \label{eq:Gm_1&Gm_2_cpx}
 \Gm_1(\z)=i e_1+e_2+\z e_3, \quad \Gm_2(z)= \z(-i e_1+e_2) - e_3. 
\end{equation}
Hence $\zD_{(x,\z)}$ is spanned by 
$\tilde{\Gm}_1(\z)_x$ and $\tilde{\Gm}_2(\z)_x$. 
Therefore the Einstein-Weyl condition is equivalent to the involutive condition 
$ [\tilde{\Gm}_1,\tilde{\Gm}_2] \in \zD. $
Proposition \ref{prop:integrability_of_complex_case} would be also proved 
in this way, it is, however, rather easier to check the integrability condition 
for $\pi^*\zD$ as we did. \\

\noindent
{\bf Definite case :}
Let $X$ be a real 3-manifold and 
$([g],\nabla)$ be a definite Weyl structure, i.e. a Weyl structure on $X$ 
with positive definite $[g]$. 
In this case, we can define complex null planes on $T_\C X$. 
If we put $\cZ=\bP(N(T^*_\C X))$, then we can define the 
complex $2$-plane distribution $\zD\subset T_\C\cZ$ in the same manner 
as the complex case by using the horizontal lift defined by 
(\ref{eq:horizontal_lift_cpx}). 
The complex conjugation $T^*_\C X\rightarrow T^*_\C X$ induces a 
fixed-point-free involution $\s: \cZ\rightarrow\cZ$
which is fiber-wise antiholomorphic. 
Notice that $\zD$ satisfies $\s^*\zD=\overline{\zD}$. 
We also define a complex $3$-plane distribution 
$\zE\subset T_\C \cZ$ by $\zE=\zD\oplus V^{0,1}$ 
where $V^{0,1}\subset T_\C\cZ$ is the ${(0,1)}$-tangent vectors 
on the fiber of $\varpi:\cZ\rightarrow X$. 
Here, we also obtain $\s^*\zE=\overline{\zE}$. 

\begin{Prop}  \label{eq:L_Riem}
 Let $([g],\nabla)$ be a definite Weyl structure on a 3-manifold $X$. 
 Let $\varpi:\cZ\to X$ be the $\CP^1$-bundle and 
 $\zE$ be the distribution on $\cZ$ constructed above. 
 Then there is a unique continuous distribution $L$ 
 of real lines on $\cZ$ which satisfies 
 $L\otimes\C= \zE\cap\overline{\zE}$ on $\cZ$. 
 Moreover the projection $\varpi(C)$ of each integral curve $C$ of $L$ 
 is a geodesic. 
\end{Prop}
\begin{proof}
 If we take a real local frame $\{e^i\}$, then we can describe 
 the situations in the similar way form (\ref{eq:conn_form_cpx}) 
 to (\ref{eq:Gm_1&Gm_2_cpx}). 
 Then $\zD=\Span\<\tilde{\Gm}_1,\tilde{\Gm}_2\>$ 
 and $\zE=\Span\<\tilde{\Gm}_1,\tilde{\Gm}_2,\pd{\bar{\z}}\>$. 
 Since $\zE+\overline{\zE}=T_\C\cZ$, $L$ exists uniquely 
 from the relation of the dimensions. 
 
 Now let us define 
 $$ l=\bar{\z}\Gm_1+\Gm_2
   = 2(\Imag \z)e_1 + 2 (\Real \z) e_2 + (|\z|^2-1)e_3. $$
 Notice that $l$ is real. 
 We can take a unique function $\g$ on $\cZ$ so that 
 $$ l^\dagger := \bar{\z}\tilde{\Gm}_1 + \tilde{\Gm}_2 +\g \pd{\bar{\z}} $$
 is real. Then we obtain $L=\Span\<l^\dagger\>$. 
 Let $p:\zE\to\zD$ be the natural projection, then 
 $p(L)=\Span\langle\,\tilde{l}\,\rangle$ 
 where $\tilde{l}=\bar{\z}\tilde{\Gm}_1 + \tilde{\Gm}_2 $. 
 By the construction, the image of an integral curve of $p(L)$ by $\varpi$ is a 
 geodesic. Pulling back to $\zE$ by $p$, we obtain the statement. 
\end{proof}

\begin{Prop} \label{prop:integrability_of_Riemannian_case}
 Let $X$ be a real 3-manifold, and $([g],\nabla)$ be a 
 definite Weyl structure on $X$ with torsion-free $\nabla$. 
 Then $([g],\nabla)$ is Einstein-Weyl if and only if 
 $\zE$ is integrable, i.e. involutive. 
\end{Prop}
\begin{proof}
 $\zE$ is integrable if and only if $\pi^*\zE$ is integrable where 
 $\pi:\cN=N(T^*_\C X)\setminus 0_X \rightarrow \cZ$. 
 If we take an orthonormal frame field $\{e_1,e_2,e_3\}$ of $T_\C X$, 
 and if we use the complex fiber coordinate $\{\l_i\}$ for $T^*_\C X$, 
 then we can define 1-forms $\c$, $\c_i$, $\t_{ij}$ on $\cN$ by 
 (\ref{eq:c&c_i&t_ij}). 
 In this case, we obtain $\pi^*\zE=\pi^*\zD+\pi^*V^{0,1}$, and 
 $\pi^*\zE = \{ v\in T^*\cN \, |\, \c(v)=\t_{ij}(v)=0 \, (\forall i,j)\}$. 
 Hence $\zE$ is integrable if and only if $\<\c,\t_{ij}\>$ is involutive. 
 By the similar arguments, this occurs if and only if $([g],\nabla)$ is Einstein-Weyl. 
\end{proof}

\begin{Rem}
 Locally speaking, $\zE / L$ defines an almost complex structure on the space of 
 geodesics on $X$. 
 Proposition \ref{prop:integrability_of_Riemannian_case} means that
 this almost complex structure is integrable if and only if 
 $([g],\nabla)$ is Einstein-Weyl (cf.\cite{bib:PT93}). 
\end{Rem}

\vspace{3mm}
\noindent
{\bf Indefinite case :}
Let $X$ be a real 3-manifold and  
$([g],\nabla)$ be a Weyl structure on $X$ whose conformal structure $[g]$ 
has signature $(-++)$. 
Let $\{e_1,e_2,e_3\}$ be a local frame field on $TX$ such that 
\begin{equation} \label{eq:-++}
 (g_{ij})=(g(e_i,e_j))=
 \begin{pmatrix}
  -1 && \\ &1& \\ &&1
 \end{pmatrix}. 
\end{equation} 
A non zero tangent vector $v\in TX$ is called 
time-like, space-like or null when $g(v,v)$ is negative, positive, or zero respectively. 
The following properties are easily checked. 
\begin{Lem} \label{Lem:space-like_and_time-like-vectors}
 \begin{enumerate}
  \item For each space-like vector $v$, there are just two real null planes 
   which contain $v$. 
  \item Each time-like vector is transverse to every real null plane. 
 \end{enumerate}
\end{Lem}

Similar to the definite case, 
we define $N(T_\C^*X)$, the space of complex null cotangent vectors, 
and $\cZ=\bP(N(T_\C^*X))$, the space of complex null planes. 
In indefinite case, we can also define $N(T^*X)$, 
the space of {\it real} null cotangent vectors, 
and $\cZ_\R=\bP(N(T^*X))$, the space of {\it real} null planes. 
There is a natural embedding $\cZ_\R\hookrightarrow\cZ$. 
The complex conjugation $T^*_\C X\rightarrow T^*_\C X$ 
induces an involution $\s: \cZ\rightarrow\cZ$ 
which is fiber-wise antiholomorphic and 
whose fixed point set coincides with $\cZ_\R$. 

Let us describe the situation explicitly using the above frame $\{e_i\}$ 
and its dual $\{e^i\}$. 
From the compatibility condition (\ref{eq:compatibility_of_Weyl_str}), 
the connection form $\w$ of $\nabla$ is written:
\begin{equation} \label{eq:conn_form}
 \w=\begin{pmatrix} 
      \p & \y^1_2 & \y^1_3 \\[1mm]
      \y^1_2 & \p & \y^2_3 \\[1mm]
	  \y^1_3 & -\y^2_3 & \p \end{pmatrix}. 
\end{equation}
We can write 
\begin{equation} \label{eq:description_of_null_cone_and_cZ}
 \begin{aligned}
   N(T^*_\C X)|_U &=\left\{ \left. 
    \sum \l_i e^i \, \right| \, -\l_1^2+\l_2^2+\l_3^2=0 \right\}, \\
  \cZ|_U &=\left\{ \left. [\l_1:\l_2:\l_3] \, \right| \, 
    -\l_1^2+\l_2^2+\l_3^2=0 \right\}. 
 \end{aligned} 
\end{equation}

 Let $U \times \CP^1\overset\sim\to \cZ_\R|_U$ be a trivialization over 
 an open set $U\subset X$ such that 
 \begin{equation} \label{eq:local_trivialization_of_cZ}
  (x,\z) \longmapsto \left[(1+\z^2)e^1+(1-\z^2)e^2+2\z e^3 \right]. 
 \end{equation}
  Here $\cZ_\R$ corresponds to 
 $\{(x,\z) \in U\times\CP^1 \, | \, \z\in\R\cup\{\infty\}\}$. 
The horizontal lift 
$\tilde{v}$ of $v\in T_xU$ at $(x,\z)\in\cZ_\R|_U$ is 
\begin{equation} \label{horizontal_lift_of_pseudo_Riemm_case}
 \tilde{v}= v + \left\{ \frac{\y^2_3+\y^1_3}{2} - \z\y^1_2 
    + \z^2 \frac{\y^2_3-\y^1_3}{2} \right\}(v) \pd{\z}. 
\end{equation}
If we define 
\begin{equation} \label{eq:Gm_1&Gm_2}
 \Gm_1(\z)=-e_1+e_2+\z e_3, \quad \Gm_2(\z)= \z(e_1+e_2) - e_3,  
\end{equation}
then $\Gm_1(\z)$ and $\Gm_2(\z)$ span the null plane 
corresponding to $(x,\z)\in\cZ_\R$. 
Define the real 2-plane distribution $\zD_\R\subset T\cZ_\R$ 
so that $\zD_\R=\Span\<\tilde\Gm_1,\tilde\Gm_2\>$ where 
$\tilde\Gm_i$ are the vector fields on $\cZ_\R$ such that 
$\tilde\Gm_{i\,(x,\z)}$ is the horizontal lift of $\Gm(\z)_x$. 

We can extend $\tilde\Gm_i$ meromorphicaly on $\cZ$, and define 
the complex 2-plane distribution $\zD\subset T_\C\cZ$ by 
$\zD=\Span\<\tilde\Gm_1,\tilde\Gm_2\>$. 
We also define a complex 3-plane distribution $\zE$ by 
$\zE=\zD\oplus V^{0,1}$ where 
$V^{0,1} \subset T_\C\cZ$ is $(0,1)$-tangent vectors.  
Then we obtain 
$$  \s^*\zD=\overline{\zD}, \quad \s^*\zE=\overline{\zE}. $$
$$ \zD_\R\otimes\C= \zD | _{\cZ_\R},\quad 
\zD_\R = \zD \cap T\cZ_\R= \zE \cap T\cZ_\R. $$

\begin{Prop} \label{prop:L_p_Riem}
 Let $([g],\nabla)$ be an indefinite Weyl structure on a 3-manifold $X$. 
 Let $\varpi:\cZ\to X$ be the $\CP^1$-bundle and $\zE$ be the distribution on $\cZ$ 
 constructed above. 
 Then there is a unique continuous distribution $L$ 
 of real lines on $\cZ$ which satisfies 
 $L\otimes\C= \zE\cap\overline{\zE}$ on $\cZ\setminus\cZ_\R$ 
 and $L\subset\zD_\R$ on $\cZ_\R$. 
 Moreover each integral curve $C$ of $L$ is contained in 
 either $\cZ\setminus\cZ_\R$ or $\cZ_\R$, and the projection 
 $\varpi(C)$ is time-like geodesic if $C\subset \cZ\setminus\cZ_\R$, 
 and null-geodesic if $C\subset\cZ_\R$. 
\end{Prop}
\begin{proof}
 Let us define a real vector field $l$ on $X$ by 
 \begin{equation}
  l=\Gm_1-\bar{\z}\Gm_2=-(1+|\z|^2)e_1+(1-|\z|^2)e_2+(\z+\bar{\z})e_3. 
 \end{equation}
 Notice that $l$ is time-like if $\Imag\z\neq 0$, 
 and null if $\Imag\z=0$. 
 We can take a unique function $\g$ on $\cZ$ so that 
 $$ l^\dagger=\tilde{\Gm}_1-\bar{\z}\tilde{\Gm_2}+\g \pd{\bar{\z}}$$
 is real. 
 Since $\tilde{l}=\tilde\Gm_1-\bar\z\tilde\Gm_2$ is real on $\cZ_\R$, 
 $\g=0$ and $l^\dagger=\tilde{l}$ on $\cZ_\R$. 
 If we put $L=\<l^\dagger\>$, then we obtain 
 $L\otimes\C=\zE\cap\overline{\zE}$ on $\cZ\setminus\cZ_\R$ 
 and $L\subset \zD_\R$ on $\cZ_\R$. 
 $L$ is unique since $E+\bar{E}=T_\C\cZ$ on $\cZ\setminus\cZ_\R$. 
 The rest statements are proved in the similar way as the definite case 
 (Proposition \ref{eq:L_Riem}). 
\end{proof}

\begin{Prop} \label{prop:integrability_of_p-Riemannian_case}
 Let $X$ be a real 3-manifold, and $([g],\nabla)$ be an 
 indefinite Weyl structure on $X$ 
 with torsion-free $\nabla$. 
 Then the following conditions are equivalent: 
 \begin{itemize}
  \item $([g],\nabla)$ is Einstein-Weyl, 
  \item the real distribution $\zD_\R$ is integrable, 
  \item the complex distribution $\zE$ is integrable. 
 \end{itemize}
\end{Prop}
\begin{proof}
If we put  
 $$ \Upsilon = -\l_1\pd{\l_1}+\l_2\pd{\l_2}+\l_3\pd{\l_3}, $$
 $$ \t_{12}=\l_1\c_2+\l_2\c_1, \quad 
  \t_{13}=\l_1\c_3+\l_3\c_1, \quad 
  \t_{23}=\l_2\c_3-\l_3\c_2 $$ 
 instead of (\ref{eq:Upsilon}) and (\ref{eq:c&c_i&t_ij}), 
 then the situation is parallel to the complex or definite case. 
\end{proof}

A direct calculation shows 
\begin{equation} \label{eq:t_23_of_indefinite_case}
\t_{23}= \l_2d\l_3- \l_3d\l_2 
  - \l_1 \(\l_1\y^2_3 + \l_2\y^1_3 - \l_3\y^1_2  \). 
\end{equation}
(\ref{eq:t_23_of_indefinite_case}) 
will be used in Section \ref{Section:Construction_of_E-W_spaces}. 

\begin{Rem}
 We can write $\zD=\<\tilde{\Gm}_1\>\oplus\<\tilde{\Gm}_2\>$ locally, 
 hence $c_1(\zD)=c_1(\<\tilde{\Gm}_1\>)+c_1(\<\tilde{\Gm}_2\>)=-2$ 
 along each $\CP^1$-fiber of $\varpi:\cZ\to X$. 
 Since $c_1(V^{0,1})=-2$, we also obtain $c_1(\zE)=-4$ along each fiber. 
\end{Rem}

\section{Hitchin correspondence}
\label{Section:Hitchin_correspondence}

In this Section, we recall the twistor correspondence for complex Einstein-Weyl structures 
introduced by Hitchin \cite{bib:Hitchin82}. 

Let $Z$ be a complex 2-manifold and $Y$ be a non-singular 
rational curve on $Z$ with the normal bundle $N_{Y/Z}\cong \cO(2)$. 
Let $X$ be the space of twistor lines, i.e. the rational curves 
which are obtained by small deformation of $Y$ in $Z$. 
By Kodaira's theorem, $X$ has a natural structure of 3-dimensional 
complex manifold, and its tangent space at $x\in X$ is identified with 
the space of sections of normal bundle $N_{Y_x/Z}$ where 
$Y_x$ is the twistor line corresponding to $x$. 

\begin{Prop} \label{prop:Hitchin_correspondence}
 There is a unique Einstein-Weyl structure on $X$ such that 
 \begin{itemize}
  \item  each non-null geodesic on $X$ corresponds to a one-parameter family 
   of twistor lines on $Z$ passing through fixed two points, and 
  \item  each null geodesic on $X$ corresponds to 
   a one-parameter family of twistor lines each of which passes through a fixed point 
   and tangents to a fixed non-zero vector there.   
 \end{itemize}
\end{Prop}
\begin{proof}
 We have $N_{Y_x/Z}\cong \cO(2)$ for each $x\in X$
 since $Y_x$ is a small deformation of $Y$. 
 We have $T_xX\cong \G(Y_x,N_{Y_x/Z})$ by definition. 
 Each holomorphic section of $N_{Y_x/Z}\simeq\cO(2)$ corresponds to  
 a degree-two polynomial $s(\z)=a\z^2+b\z+c$ 
 where $\z$ is the inhomogeneous coordinate on $Y_x$. 
 We can define the conformal structure $[g]$ so that 
 a tangent vector in $T_xX$ is null if and only if the corresponding 
 polynomial $s(\z)$ has double roots, i.e. when $b^2-4ac=0$. 

 If we fix, maybe infinitely near, two points in $Z$, 
 then the twistor lines passing through these points make  
 a one-parameter family. 
 This family corresponds to a 
 holomorphic curve on $X$. 
 Let $\cF$ be the family of such holomorphic curves. 
 Then, by Proposition \ref{prop:family_of_curves=>proj_str}, 
 we obtain unique projective structure $[\nabla]$ on $X$ 
 such that $\cF$ coincides to the geodesics. 

 Now, we prove that there is a unique torsion-free $\nabla\in[\nabla]$ such that 
 $([g],\nabla)$ defines a Weyl structure. 
 For this purpose, we first fix an arbitrary torsion-free $\nabla\in[\nabla]$, and 
 check that the second fundamental form on 
 each null surface with respect to $\nabla$ vanishes. 

 For each point $p\in Z$, the two-parameter family of twistor lines 
 passing through $p$ corresponds to a null surface $S$ on $X$. 
 Notice that $S$ is totally geodesic and 
 naturally foliated by null geodesics each of which corresponds to 
 a tangent line at $p$. 
 Let $N=TX|_S/TS$ be the normal bundle of $S$.
 The second fundamental form $I\! I : TS\otimes TS \to N$ is defined by 
 $v\otimes w \to [\nabla_v w]^N$ 
 where the value does not depend on how to extend $w$. 
 Take a frame field $\{e_1,e_2,e_3\}$ on $TX|_S$ so that $e_1$ is null and 
 $TS=\<e_1,e_2\>$.  Then the metric tensor is 
 $$ g=(g_{ij})=\begin{pmatrix} 0&0&* \\ 0 &*&* \\ *&*&*
   \end{pmatrix}. $$ 
 Since $\nabla$ is torsion-free, 
 $\nabla_{e_1}e_2-\nabla_{e_2}e_1=[e_1,e_2] \in TS$, 
 so $g(\nabla_{e_1}e_2,e_1)=g(\nabla_{e_2}e_1,e_1)$. 
 Since $g_{13}\neq 0$, we obtain
 \begin{equation} \label{eq:G^3_{12}=G^3_{21}}
   \G^3_{12}=\G^3_{21}. 
 \end{equation}
 On the other hand $S$ is totally geodesic, we obtain 
 $$ 0=g(\nabla_\x \x, e_1)= \x^1\x^2g_{13}\( \G^3_{12}+\G^3_{21}\), $$
 for every tangent vector  $\x=\x^1e_1+\x^2e_2$ on $S$. 
So we obtain $\G^3_{12}+\G^3_{21} =0$, and combining with 
(\ref{eq:G^3_{12}=G^3_{21}}), we obtain $\G^3_{12}=\G^3_{21}=0$. 
 Hence $g(\nabla_\x\y,e_1)=0$ for every vector field 
 $\x$ and $\y$ on $S$, and this means $I\! I=0$ on $S$. 

 Next we claim that there are functions $a_i,b_i \ (i=1,2,3)$ on $X$ such that 
 \begin{equation} \label{eq:nabla_g_ijk}
 (\nabla g)_{ijk}=a_i g_{jk}+\frac{1}{2}b_j g_{ik} + \frac{1}{2} b_k g_{ij}. 
 \end{equation}
 Since $I\! I=0$ for every null surface, we obtain 
 \begin{equation} \label{eq:key_eq_from_II=0}
  \nabla_\y g(\x,\x)=0
 \end{equation}
 for every null vector $\x$ and every vector $\y$ satisfying $g(\y,\x)=0$. 
 Let us fix a local frame $\{e_i\}$ on $X$. 
 If we put $\x=\x^i e_i, \y=\y^i e_i \ (i=1,2,3)$ and 
 $\vp_{ijk}=\nabla_{e_i}(e_j,e_k)$, then 
 (\ref{eq:key_eq_from_II=0}) is written 
 \begin{equation} \label{eq:vanishing_of_the_contraction}
  (\vp_{ijk}\x^j\x^k)\y^i=0. 
 \end{equation}
 Since $\x$ runs all null vectors, $(\x^i)$ moves the conic 
 $$ C= \left\{ \left. [\x^1:\x^2:\x^3]\in \CP^2 \ 
   \right| \ \x^i\x^jg_{ij}=0 \right\}. $$
 For fixed $\x$, $(\y^i)$ moves the line 
 $$ L(\x)=\{ \left. [\y^1:\y^2:\y^3]\in \CP^2 \ 
   \right| \ \y^i(\x^jg_{ij})=0  \}. $$
 Since (\ref{eq:vanishing_of_the_contraction}) holds for every 
 $[\y^i]\in L(\x)$, we can take a function $b(\x)$ satisfying 
 $$ \vp_{ijk}\x^j\x^k= b(\x) \x^j g_{ij} $$ 
 for every $\x\in C$ and $i=1,2,3$. 
 Then we can take $b(\x)$ to be a degree-one polynomial. 
 Actually, since $\x^jg_{ij}\ (i=1,2,3)$ does not vanish at once, 
 $b(\x)=(\vp_{ijk}\x^j\x^k) / (\x^j g_{ij}) $ defines holomorphic section of 
 $\cO(1)$ over $\CP^2$. 
 If we put $b(\x)=b_k\x^k$, then we obtain 
 $$ (\vp_{ijk}-b_kg_{ij})\x^j\x^k=0 $$
 for $i=1,2,3$. Here $b_k \ (k=1,2,3)$ are functions on $X$.  
 Since these equation hold for every $\x\in C$, there are 
 functions $a_i$ on $X$ such that 
 $$ (\vp_{ijk}-b_kg_{ij})X^jX^k=a_i(g_{jk} X^j X^k) $$
 for every $(X^j)\in\C^3$ and $i=1,2,3$. 
 Noticing the symmetry, we obtain (\ref{eq:nabla_g_ijk}). 

 Finally, if we define a new connection $\tilde{\nabla}$ by 
 \begin{equation}
   \tilde{\G}^i_{jk}=\G^i_{jk}+\frac{1}{2}b_j + \frac{1}{2} b_k, 
 \end{equation}
 then $\tilde{\nabla}\in[\nabla]$ and $\tilde{\nabla}$ satisfies 
 $$ (\tilde{\nabla} g)_{ijk}=(a_i-b_i) g_{jk}, $$
 i.e. $\tilde{\nabla}$ is compatible to $[g]$. 
 Moreover, $([g],\tilde{\nabla})$ is Einstein-Weyl 
 since the integrable condition in Proposition \ref{prop:integrability_of_complex_case}
 is automatically satisfied from the construction. 
 Notice that such connection is unique since the compatibility condition is not 
 satisfied for any other torsion-free connection in $[\nabla]$. 
\end{proof}

\begin{Rem} \label{Rem:double_fibration_complex_case}
 Let $\cX=\{(x,p)\in X\times Z \, | \, p\in Y_x\}$, then we obtain the 
 double fibration $X \overset\varpi\leftarrow \cX \overset \Gf\to Z$ 
 where $\varpi$ and $\Gf$ are the projections. 
 Each $u\in \cX$ defines a null plane at $\varpi(u)\in X$ as a tangent plane of 
 the null surface corresponding to $\Gf(u)\in Z$. 
 Hence we obtain a natural map $\cX\to\cZ=\bP(N(T^{*\,1,0}_\C X))$ 
 which is in fact biholomorphic. 
 Identifying $\cX$ with $\cZ$, we obtain 
 $\zD=\ker\{\Gf_*: T^{1,0}_\C\cX \to T^{1,0}_\C Z \}$. 
\end{Rem}

Hitchin introduced two examples of Einstein-Weyl spaces each of which is obtained from 
a complex twistor space \cite{bib:Hitchin82}. 
The twistor space of one of them is 
$$ Z=\left\{ [z_0:z_1:z_2:z_3]\in\CP^3 \, | \, z_1^2+z_2^2+z_3^3=0 \right\}. $$ 
In this case, the twistor lines are the plane sections, and the corresponding 
Einstein-Weyl space is flat. 
In the other case, the twistor space is 
\begin{equation} \label{eq:standard_complex_twistor_space}
 Z=\left\{ [z_0:z_1:z_2:z_3]\in\CP^3 \, | \, z_0^2+z_1^2+z_2^2+z_3^3=0 \right\}. 
\end{equation}
In this case, the twistor lines are also the plane sections, and the corresponding 
Einstein-Weyl space is constant curvature space. 
We study more detail about the latter one in next Section. 

\section{The standard case} 
\label{Section:The_standard_case}

In this section, the standard model of LeBrun-Mason type correspondence is explained. 
We start from Hitchin's example (\ref{eq:standard_complex_twistor_space}), 
and construct the model case as a real slice of it (c.f.\cite{bib:PT93}). 

If we change the coordinate, (\ref{eq:standard_complex_twistor_space}) 
can be written $\{[z_i]\in\CP^3 \, | \, z_0z_3=z_1z_2 \}$ 
which coincides with the image of Segre embedding 
$\CP^1\times\CP^1\hookrightarrow \CP^3$ 
$$ ([u_0:u_1],[v_0:v_1]) \longmapsto 
  [u_0v_0:u_0v_1:u_1v_0:u_1v_1]. $$
So we usually denote $Z=\CP^1\times\CP^1$. 
Since the twistor lines are the plane sections, 
the twistor lines are parametrized by $X=\CP^{* 3}$. 
We introduce a homogeneous coordinate $[\x^i]\in \CP^{* 3}$ 
so that $[\x^i]$ corresponds to the plane $\{[z_i]\in\CP^3 \, | \, \x^iz_i=0\}$. 
Let 
$$ X_{\text{sing}}=
  \left\{ \left. [\x^i]\in \CP^{* 3} \, \right| \, \x^0\x^3=\x^1\x^2 \right\} $$ 
be the set of planes each of which tangents to $Z$. 
If $[\x^i]\in X_{\text{sing}}$, then the plane section degenerates to 
two lines 
$$ \(\CP^1\times [-\x^1:\x^0]\) \cup \( [-\x^2:\x^0]\times\CP^1 \) $$
intersecting at the tangent point. 
We call such a plane section a {\it singular twistor line} on $Z$. 
Since Proposition \ref{prop:Hitchin_correspondence} does not work on 
$X_{\text{sing}}$, the Einstein-Weyl structure is defined only on 
$X\setminus X_{\text{sing}}$. 

\vspace{2mm}
Next we introduce real structures, i.e. antiholomorphic involution on $Z$. 
There are several ways to introduce such structure. 
For example, if we take the fixed-point-free involution 
$$ \s' : ([u_0:u_1],[v_0:v_1]) \longmapsto 
  ([\bar{u}_1:\bar{u}_0],[\bar{v}_1:-\bar{v}_0]), $$
then $\s'$ extends to the involution on $\CP^3$ by 
$$ [z_0:z_1:z_2:z_3]\longmapsto
  [\bar{z}_3:-\bar{z}_2:-\bar{z}_1:\bar{z}_0]. $$
Then we also obtain an antiholomorphic involution on $X$, 
and let $X_\R$ be its fixed point set. 
Since $X_\R\cap X_{\text{sing}}$ is empty, 
we obtain a real Einstein-Weyl structure on whole $X_\R\cong\RP^3$ 
as a real slice of the complex Einstein-Weyl structure on $X\setminus X_{\text{sing}}$. 
This is nothing but the definite Einstein-Weyl structure 
induced from the standard constant curvature metric on $\RP^3$. 

\vspace{2mm}
Our main interest is, however, in the indefinite case. Let 
$$ \s : ([u_0:u_1],[v_0:v_1]) \longmapsto 
  ([\bar{v}_1:\bar{v}_0],[\bar{u}_1:\bar{u}_0]), $$ 
be another involution on $Z$ whose fixed point set is denoted by $Z_\R$. 
$\s$ extends to the involution on $\CP^3$ by 
$$ [z_0:z_1:z_2:z_3]\longmapsto
  [\bar{z}_3:\bar{z}_1:\bar{z}_2:\bar{z}_0]. $$
Then we also obtain an involution on $X$, 
and let $X_\R$ be its fixed point set. In this case, 
$X_{\R,\text{sing}}=X_\R\cap X_{\text{sing}}$ is nonempty. 

Let $(\y_1,\y_2)=(u_0/u_1,v_0/v_1)$ be a coordinate 
on $Z=\CP^1\times\CP^1$, 
and let us denote $\t(\y)=\bar{\y}^{-1}$. 
Then $\s(\y_1,\y_2)=(\t(\y_2),\t(\y_1))$ 
and $Z_\R=\{ (\y,\t(\y)) \, | \, \y\in\CP^1 \}. $
In this coordinate, each non-singular twistor line $l$ is written as a graph of 
some M\"obius transform $f:\CP^1\to\CP^1$, i.e. 
$l=\{(\y,f(\y)) \, | \, \y\in\CP^1 \}$. 
$l$ is $\s$-invariant iff 
$\t (f(\y))=f^{-1}(\t(\y))$, and then we can write  
$$ f(\y)=\frac{A\y-B}{\bar{B}\y-C} $$
for some $(A,B,C)\in \R\times\C\times\R$ satisfying $|B|^2-AC\neq 0$. 
$l\cap Z_\R$ is nonempty if $|B|^2-AC>0$, and is empty if $|B|^2-AC<0$. 

In the non-singular case, 
the parameters $(A,B,C)$ can be normalized so that $|B|^2-AC=\pm 1$. 
Since $(A,B,C)$ and $(-A,-B,-C)$ defines the same M\"obius transform, 
we obtain $X_\R \setminus X_{\R,\text{sing}}\cong H\sqcup H'$ where
$$ \begin{aligned}
  H &=\left\{(A,B,C)\in \R\times\C\times\R \, \left| \, 
      |B|^2-AC=1 \right. \right\}/ \pm, \\ 
  H' &=\left\{(A,B,C)\in \R\times\C\times\R \, \left| \, 
      |B|^2-AC=-1 \right. \right\}/ \pm. 
  \end{aligned} $$
We obtain an indefinite Einstein-Weyl structure on $H$ 
and a definite Einstein-Weyl structure on $H'$ 
as a real slice of $X\setminus X_{\text{sing}}$ . 
The conformal structures are the class of 
$$g=|dB|^2-dAdC$$
which is indefinite on $H$ and definite on $H'$. 

If we identify $ \CP^1\overset\sim\to Z_\R : 
\w \mapsto \(\w,\bar{\w}^{-1}\), $
then each twistor line corresponding to $[A,B,C]\in H$ intersects with $Z_\R$ 
by the circle
\begin{equation} \label{eq:standard_circle}
 \left\{ \w\in\CP^1 \, \left| \, A|\w|^2-B\bar{\w}-\bar{B}\w+C=0 
  \right. \right\}. 
\end{equation}
Hence $H$ is naturally identified with the set of circles on $\CP^1$, 
and its double cover 
$$ \tilde{H}=\left\{(A,B,C)\in \R\times\C\times\R \, \left| \, 
      |B|^2-AC=1 \right. \right\} \cong S^2\times \R$$ 
is identified with the set of oriented circles on $\CP^1$.  
Since each circle divides the twistor line into two holomorphic disks, 
$\tilde{H}$ is identified with the set of holomorphic disks in $Z$ 
whose boundaries lie on $Z_\R$. 

There is a natural action of $\PSL(2,\C)$ on $H$, $H'$ and $\tilde{H}$
defined by the following way. 
Each $\p\in\PSL(2,\C)=\Aut(\CP^1)$ induces an automorphism on $Z$ by 
\begin{equation} \label{eq:PSL(2,C)_action}
 \p_*:(\y_1,\y_2)\mapsto (\p(\y_1),\t\p\t(\y_2)). 
\end{equation}
$\p_*$ maps each $\s$-invariant twistor line 
to the other $\s$-invariant twistor line. 
Since $\p_*$ preserves $Z_\R$, $\p_*$ preserves $H$ and $H'$. 
Obviously this action lifts to an automorphism on $\tilde{H}$, 
and we will see later that this action on $\tilde{H}$ is transitive. 

\vspace{2mm}
Now, we introduce an explicit description of the holomorphic disks 
corresponding to $\tilde{H}$. 
Let $M\cong \CP^1\times \R = U_1\cup U_2$ where
$U_i=\{(\l_i,t)\in \C\times \R\}$ are patched by $\l_2=\l_1^{-1}$. 
Let $\varpi:\cX_+\to M$ be the disk bundle  
$$ \begin{aligned}
 {} & \cX_+ =  (U_1\times \bD) \cup (U_2\times \bD), \\
 {} & (\l_1,t;z_1) \sim (\l_2,t;z_2) \, \Longleftrightarrow \,
  \l_2=\l_1^{-1},\, z_2=\frac{\bar{\l}_1}{\l_1}z_1, 
\end{aligned} $$ 
where $\bD=\{ z\in\C\, |\, |z|\leq 1\}$. 
We denote $\cX_\R=(U_1\times\del \bD)\cup(U_2\times\del \bD)$, 
and notice that $\cX_\R$ is a circle bundle with $c_1(\cX_\R)=2$ along each 
fiber of $\varpi$. 
Let us define a smooth map $\Gf: \cX_+ \to Z$ by 
$$ \begin{aligned}
  U_1\times \bD \ni (\l_1,t;z_1) & \longmapsto 
   \( \frac{z_1+r\l_1}{-\bar{\l}_1z_1+r}, \frac{rz_1-\l_1}{r\bar{\l}_1z_1+1} \), \\
  U_2\times \bD \ni (\l_2,t;z_2) & \longmapsto 
   \( \frac{\bar{\l}_2z_2+r}{-z_2+r\l_2}, \frac{r\bar{\l}_2z_2-\l}{rz_2+\l_2} \), 
 \end{aligned} $$
where $r=e^t$. 
In this way, we have obtained the following double fibration; 
\begin{equation} \label{eq:standard_double_fibration}
  \xymatrix{ & \cX_+ \ar[dl]_\varpi \ar[dr]^\Gf & \\ 
     M && Z}
\end{equation}
We use the coordinate $\l\in\C\cup \{\infty\}= \CP^1$ satisfying  
$\l=\l_1$ on $U_1$, and we define $D_{(\l,t)}=\Gf\action\varpi^{-1}(\l,t)$. 
Then $\{D_{(\l,t)}\}_{(\l,t)\in M}$ gives the family of holomorphic disks 
which coincides with the family corresponding to $\tilde{H}$ above. 
Hence naturally $M\cong \tilde{H}$. 
Notice that we arranged so that the center of $D_{(\l,t)}$, i.e. $z=0$, lies on 
$$Q=\left\{ (\l,-\l) \in Z \, \left| \, \l\in\CP^1 \right. \right\}$$
which is a twistor line on $Z$ corresponding to $[1,0,1]\in H'$. 

\vspace{2mm}
We have already defined a $\PSL(2,\C)$-action on $M=\tilde{H}$ 
by (\ref{eq:PSL(2,C)_action}). 
For each element $\p\in\PSU(2)\subset\PSL(2,\C)$, we can check that 
$\p_*(D_{(\l,t)})=D_{(\p(\l),t)}$. 
Since $\PSU(2)$ acts transitively on $\CP^1$, 
$\PSU(2)$ acts transitively on $\CP^1\times\{t\}\subset M$
for each $t\in\R$. 
On the other hand,  
\begin{equation} \label{eq:action_changing_t}
 \p=\left[ \begin{array}{cc} e^{-t} & \\ & e^t \end{array} \right] 
  \in \PSL(2,\C) 
\end{equation}
gives the automorphism $\p_*$ which maps 
the disk $D_{(0,1)}$ to $D_{(0,2t)}$. 
Hence the action of $\PSL(2,\C)$ on $M=\tilde{H}$ is transitive. 

Let $S(TZ_\R)=(TZ_\R\setminus 0_{\Z_\R})/\R_+$ be the circle bundle on $\Z_\R$ 
where $0_{\Z_\R}$ is the zero section and $\R_+$ is positive real numbers 
acting on $TZ_\R$ by scalar multiplication.  
On $\cX_\R$, we can take a nowhere vanishing vertical vector field {\bf v}, 
i.e. $\varpi_*(\text{\bf v})=0$, so that the orientation matches to 
the complex orientation of the fiber of $\varpi:\cX_+\to M$. 
Since $\Gf_*(\text{\bf v})$ does not vanish anywhere, we can define a smooth map 
$\tilde{\Gf}:\cX_\R\to S(TZ_\R)$ by $u\mapsto [\Gf_*({\bf v}_u)]$. 
Then we obtain the following diagram: 
$$ \xymatrix{ 
   \cX_\R \ar[rr]^(0.45){\tilde{\Gf}} \ar[dr]_\Gf && S(TZ_\R) \ar[dl] \\
   & Z_\R &} $$

\begin{Prop} \label{Prop:properties_over_S_t}
 Let $S_t=\CP^1\times \{t\} \subset M$, 
 and let $\Gf_t$ and $\tilde{\Gf}_t$ be the restriction of $\Gf$ and $\tilde{\Gf}$ on 
 $\varpi^{-1}(S_t)$ respectively. Then, for each $t\in \R$, 
 \begin{enumerate}
  \item $\Gf_t:(\cX_+\setminus\cX_\R)|_{S_t}\to Z\setminus Z_\R$  
    is diffeomorphic, 
  \item $\tilde{\Gf}_t:\cX_\R|_{S_t}\to S(TZ_\R)$ is diffeomorphic, 
  \item $\Gf_t:\cX_\R|_{S_t}\to Z_\R$ is an $S^1$-fibration such that 
    each fiber is transverse to the vertical distribution of $\varpi: \cX_\R\to M$, 
 \end{enumerate}
 In particular, $\{D_{(\l,t)}\}_{\l\in\CP^1}$ 
 gives a foliation on $Z\setminus Z_\R$ for each $t\in\R$. 
\end{Prop}

\begin{Rem}
 Notice that, from {\it 2} above, the following holds:  
 for each $t\in\R$, $p\in Z_\R$ and non zero $v\in T_p Z_\R$, there is a unique 
 $x\in S_t$ such that $p\in\del D_x$ and $v\parallel D_x$ 
 (cf. Definition \ref{Def:adapted_vector}). 
\end{Rem}

\begin{proof}[Proof of \ref{Prop:properties_over_S_t}]
 We can assume $t=0$ by changing the parameter $t\in\R$
 by the automorphism of type (\ref{eq:action_changing_t}). 

 When $t=0$, we can interpret the situation to a geometry on $S^2$
 in the following way. 
 Let $S^2=\{(x_1,x_2,x_3)\in\R^3 \,|\, \sum x_i^2=1\}$ and 
 $\Gp : \CP^1\overset\sim \to S^2$ be the stereographic projection, 
 $$\Gp : \l \longmapsto \( \frac{2\Real \l}{1+|\l|^2}, \frac{2\Imag \l}{1+|\l|^2}, 
   \frac{1-|\l|^2}{1+|\l|^2}, \). $$
 We identify $Z$ with $S^2\times S^2$ by the diffeomorphism 
 $Z\overset\sim \to S^2\times S^2 : 
   (\y_1,\y_2)\mapsto (\Gp(\y_1),\Gp\action\t(\y_2))$ 
 where $\t(\y)=\bar{\y}^{-1}$.    
 Notice that $Z_\R$ corresponds to the diagonal in this identification. 

 Recall that $D_{(\l,0)}$ is the image of $\bD\to Z :$ 
 $$ z\longmapsto (\y_1,\y_2)=
 \( \frac{z+\l}{-\bar{\l}z+1}, \frac{z-\l}{\bar{\l}z+1}\). $$ 
 Then $\del D_{(\l,0)}\subset Z_\R$ corresponds to the 
 big circle on the diagonal $S^2\subset S^2\times S^2$ cut out by the plane 
 \begin{equation}  \label{eq:big_circle}
   2(\Real \l) x_1 + 2(\Imag \l) x_2 + (1-|\l|^2) x_3 =0. 
 \end{equation}
 Hence we obtain the one-to-one correspondence between 
 $\l\in\CP^1$ and the oriented big circle $\Gp(\del D_{(\l,0)})$, 
 where the orientation is induced from the natural orientation of $\Gp(D_{(\l,0)})$. 
 Moreover, we claim that the following conditions are equivalent: 
 \begin{list}{(A\,\arabic{mynum})}%
  {\usecounter{mynum} \itemsep 0in  \leftmargin .4in}
   \item $(\y_1,\y_2)\in Z$ lies on $D_{(\l,0)}$, 
   \item the oriented big circle $\Gp(\del D_{(\l,0)})$ rounds anti-clockwise around 
     $\Gp(\y_1)$, and this big circle coincides with the set of points on $S^2$ 
	 which have the same distant from 
     $\Gp(\y_1)$ and $\Gp\action\t(\y_2)$ with respect to the standard 
     Riemannian metric on $S^2$. 
 \end{list}
 Indeed, if $(\y_1,\y_2)\in D_{(\l,0)}$, then the point 
 $$\Gp(\y_1) + \Gp\action\t(\y_2) \in \R^3$$ 
 lies on the plane (\ref{eq:big_circle}), 
 hence the big circle $\Gp(\del D_{(\l,0)})$ satisfies (A2). 
 The converse is easy. 
 In particular, the following conditions are equivalent : 
 \begin{list}{(B\,\arabic{mynum})}%
  {\usecounter{mynum} \itemsep 0in  \leftmargin .4in}
    \item $(\y_1,\y_2)\in Z_\R$ lies on $\del D_{(\l,0)}$, 
    \item the big circle $\Gp(\del D_{(\l,0)})$ passes through 
    $\Gp(\y_1)=\Gp\action\t(\y_2)$. 
 \end{list}
 The statement follows directly from these interpretation. 
 Actually, for each $p=(\y_1,\y_2)\in Z\setminus Z_\R$, 
 the big circle satisfying (A2) exists uniquely, hence {\it 1} holds. 
 For each $p=(\y_1,\y_2)\in Z_\R$, 
 $S(T_pZ_\R)$ corresponds to the oriented big circles satisfying (B2), 
 hence {\it 2} and {\it 3} follows. 
\end{proof}

The geometry on $M$ is characterized by the double fibration 
(\ref{eq:standard_double_fibration}) in the following way. 
\begin{Prop} \label{prop:geometry_of_the_standard_double_fibration}
 \begin{enumerate}
  \item For each $p\in Z_\R$, 
   $\GS_p=\{ x\in M \, | \, p\in \del D_x \}=\varpi\action\Gf^{-1}(p)$ 
   is maximal connected null surface on $M$ 
   and every null surface can be written in this form. 
  \item For each $p\in Z\setminus Z_\R$, 
   $\GC_p=\{ x\in M \, | \, p\in D_x \}=\varpi\action\Gf^{-1}(p)$ 
   is maximal connected time-like geodesic on $M$ 
   and every time-like geodesic can be written in this form. 
  \item For each $p\in Z_\R$ and each non-zero $v\in T_p Z_\R$, 
   $\GC_{p,v}=\{x\in M \, | \, p\in \del D_x, v\parallel D_x \}
     =\varpi\action\tilde\Gf^{-1}([v])$  
   is maximal connected null geodesic on $M$ 
   and every null geodesic can be written in this form. 
  \item For each distinguished $p, q\in Z_\R$, 
   $\GC_{p,q}=\{x\in M \, | \, p,q\in \del D_x \}=\GS_p\cap \GS_q$ 
    is closed connected space-like geodesic on $M$ 
   and every space-like geodesic can be written in this form. 
 \end{enumerate}
\end{Prop}
\begin{proof}
 Since $\{\del D_{(\l,t)}\}$ is the set of oriented circles of the form 
 (\ref{eq:standard_circle}), we obtain 
 \begin{itemize}
  \item $\GS_p\simeq S^1\times\R$ for each $p\in Z_\R$,  
  \item $\GC_{p,v}\simeq \R$ for each $p\in Z_\R$ and 
     non zero vector $v\in T_p Z_\R$. 
  \item $\GC_{p,q}\simeq S^1$ for each distinguished $p,q\in Z_\R$, 
 \end{itemize}
 Since $\GS_p$ is a real slice of complex null surface, it is real null surface. 
 Moreover, it is a maximal connected null surface since $\GS_p$ is closed in $M$. 
 Hence {\it 1} holds. 
 In the similar way, we can see that $\GC_{p,v}$ is a maximal connected real 
 null geodesic, so {\it 3} holds. 
 $\GC_{p,q}$ is also a maximal connected real non-null geodesic. 
 Notice that $\GC_{p,q}$ is contained in the null surface $\GS_p$. 
 Since null plane never contain time-like vector, $\GC_{p,q}$ is a space-like 
 geodesic (cf. Lemma \ref{Lem:space-like_and_time-like-vectors}). 
 Hence {\it 4} holds. 

 Now we check {\it 2}. Let $p\in Z\setminus Z_\R$ and notice that 
 every $\s$-invariant twistor line passing through $p$ also passes through $\s(p)$. 
 So $\GC_p$ is a real slice of the complex geodesic corresponding to the two points 
 $p,\s(p)$. Hence $\GC_p$ is a geodesic. 
 From Proposition \ref{Prop:properties_over_S_t}, we obtain $\GC_p\simeq\R$ 
 which is closed in $M$. 
 Hence $\GC_p$ is maximal connected geodesic. 
 To see that $\GC_p$ is a time-like geodesic, 
 it is enough to check that $\GC_p$ is transversal to every null plane at each point 
 (cf. Lemma \ref{Lem:space-like_and_time-like-vectors}). 
 Notice that, if we fix three points on $Z$, 
 there are at most one twistor line containing them, 
 hence $\GC_p\cap\GS_q=\{x\in M \, |\, p,\s(p),q\in D_x\}$ is at most one point 
 for each $q\in Z_\R$. Thus $\GC_p$ is time-like. 
\end{proof}

In particular, we obtain the following. 
\begin{Cor}
 The indefinite Einstein-Weyl structure on $M$ 
 constructed above is space-like Zoll. 
\end{Cor}

Let $\cX=\cX_+\cup_{\cX_\R} \cX_-$ be a $\CP^1$ bundle over $M$
where $\cX_-=\overline{\cX_+}$ is the copy of 
$\cX_+$ with fiber-wise opposite complex structure. 
On the other hand, we have a $\CP^1$-bundle $\cZ$ on $M$ 
equipped with the distributions $\zD_\R,\zE,L$ and so on. 
Then, similar to Remark \ref{Rem:double_fibration_complex_case}, 
there is a natural identification $\cX\overset\sim\to\cZ$
such that 
\begin{itemize}
 \item for each $p\in Z_\R$, $\Gf^{-1}(p)$ corresponds to 
      an integral surface of $\zD_\R$, 
 \item for each $p\in Z\setminus Z_\R$, $\Gf^{-1}(p)$ corresponds to 
      an integral curve of $L$ in $\cX\setminus\cX_\R$, 
 \item for each $p\in Z_\R$ and $[v]\in S(T_pZ_\R)$, $\tilde\Gf^{-1}([v])$ 
      corresponds to an integral curve of $L$ in $\cX_\R$. 
\end{itemize}
Hence the following holds: 
\begin{itemize}
 \item $\zD_\R=\zE\cap T\cX_\R=\ker\{\Gf_*: T\cX_\R \to TZ_\R\}$ 
   on $\cX_\R$, 
 \item $L=\ker\{\Gf_*: T\cX \to TZ\}$ on $\cX_+\setminus\cX_\R$, and 
 \item $L=\ker\{\tilde\Gf_* : T\cX_\R\to S(TZ_\R)\}$ on $\cX_\R$.
\end{itemize}

Recall that we denote $S_t=\CP^1\times\{t\}$, and let us denote 
$\cX_t=\varpi^{-1}(S_t)$ where $\varpi:\cX\to M$ is the projection. 
Let $\zE_t=\zE\cap T_\C\cX_t$ for each $t$. 
Then, since $L\cap T\cX_t=0$, we obtain $\zE=(L\otimes\C)\oplus\zE_t$.  
From $\zE\cap\overline{\zE}=L\otimes\C$ and 
$\zE\oplus\overline{\zE}=T\cX$, we obtain 
$\zE_t\oplus\overline{\zE_t}=T\cX_t$. 
Moreover, since $\zE$ is integrable, $\zE_t$ is also integrable. 
Hence $\zE_t$ defines a complex structure on $\cX_t$. 

Now we claim that $\Gf_t: (\cX_+\setminus\cX_\R)|_{S_t}\to Z\setminus\Z_\R$ 
is holomorphic with respect to the above complex structure. 
Consider the complex Einstein-Weyl space
$M_\C=X\setminus X_{\text{sing}}$ defined at the beginning of this section, 
and let $\cZ_\C=\bP(N(T^{*\, 1,0}M_\C))$. 
Then we obtain the double fibration $M_\C \leftarrow \cZ_\C \overset{\Gf_\C}\to Z$ 
where $\Gf_\C$ is holomorphic. 
On the other hand, there is natural embedding 
$i_t: (\cX_+\setminus \cX_\R)|_{S_t} \to \cZ_\C$ 
which is holomorphic since it preserves the distributions. 
Since $\Gf_t=\Gf_\C\action i_t$, 
$\Gf_t$ is holomorphic on $(\cX_+\setminus \cX_\R)|_{S_t}$. 

From the above argument, we obtain 
$\zE_t=(\Gf_t)^{-1}_*(T^{0,1}Z) \subset \Gf_*^{-1}(T^{0,1}Z)$ on 
$\cX_+\setminus\cX_\R$. 
Since $L\otimes\C=\ker\Gf_*$ there,  we obtain 
$\zE=(L\otimes\C)\oplus\zE_t \subset \Gf_*^{-1}(T^{0,1}Z)$ on  
$\cX_+\setminus\cX_\R$. 
Then we also have $\overline{\zE}\subset\Gf_*^{-1}(T^{1,0}Z)$. 
Since $\zE+\overline{\zE}=T_\C\cX_+$ and
$\zE\cap\overline{\zE}=L\otimes\C$, we obtain 
$\zE= \Gf_*^{-1}(T^{0,1}Z)$ on $\cX_+\setminus\cX_\R$. 

In this way, we have proved the following. 
\begin{Prop} 
Identifying $\cX=\cX_+\cup\cX_-$ with $\cZ$,  
\begin{enumerate}
 \item $\zE=\Gf_*^{-1}(T^{0,1}Z)$ on $\cX_+$ where 
   $\Gf_* : T_\C\cX_+ \to T_\C Z$, 
 \item $\zD_\R=\zE\cap T\cX_\R=\ker\{\Gf_*: T\cX_\R \to TZ_\R\}$ 
   on $\cX_\R$, 
 \item $L=\ker\{\Gf_*: T\cX_+ \to TZ\}$ on $\cX_+\setminus\cX_\R$, and 
 \item $L=\ker\{\tilde\Gf_* : T\cX_\R\to S(TZ_\R)\}$ on $\cX_\R$.
\end{enumerate}
\end{Prop}

\vspace{2mm}
It is convenient considering the compactification of $M$ and $\cX_+$. 
Let $I=[-\infty,\infty]$ be the natural compactification of $\R$. 
If we put $\hat{M}=\CP^1\times I$, 
then we obtain a natural embedding $\iota:M\hookrightarrow\hat{M}$. 
Next, let $\Psi: \cX_+ \to \hat{M} \times Z$ be the embedding defined by 
$\Psi(u) = \(\iota\action\varpi(u), \Gf(u)\)$. 
Let us define $\hat{\cX}_+$ and $\hat{\cX}_\R$ as the 
closure of $\Psi(\cX_+)$ and $\Psi(\cX_\R)$ in $\hat{M} \times Z$. 
Then we obtain the double fibration 
\begin{equation} \label{eq:compactified_double_fibration}
 \xymatrix{ & ( \hat{\cX}_+, \hat{\cX}_\R)  
 \ar[dl]_{\hat{\varpi}} \ar[dr]^{\hat{\Gf}}& \\
    \hat{M} & & (Z,Z_\R)} 
\end{equation}	
where $\hat{\varpi}$ and $\hat\Gf$ is the projections. 

Notice that $\hat{\varpi}^{-1}(x)$ is no longer a disk 
for $x=(\l,\pm\infty)\in\del\hat{M}$, but a {\it marked} $\CP^1$ 
whose marking point is $\hat{\varpi}^{-1}(x)\cap\hat{\cX}_\R$. 
We denote them by 
\begin{equation} \label{eq:marked_CP^1}
 \begin{aligned}
   D_{(\l,\infty)} &= \hat{\varpi}^{-1}(\l,-\infty) = \{\l\}\times\CP^1, \\
   D_{(\l,-\infty)} &=\ \hat{\varpi}^{-1}(\l,\infty) \ = \CP^1\times\{-\l\},
  \end{aligned} 
\end{equation}
where $D_{(\l,\infty)}$ is marked at $(\l,\bar{\l}^{-1})$ and 
$D_{(\l,-\infty)}$ is marked at $(-\bar{\l}^{-1},-\l)$. 

Recall the definitions of $\GC_p$, $\GC_{p,v}$ and so on 
introduced at Proposition \ref{prop:geometry_of_the_standard_double_fibration}. 
We define $\hat{\GC}_p$, $\hat{\GC}_{p,v}$ and so on 
as the compactification in $\hat{M}$.  
Then the following properties are easily checked. 
\begin{Prop} \label{prop:foliation_on_Z_R}
 \begin{enumerate}
  \item For each $p\in Z\setminus Z_\R$, 
    $\hat{\cX}_\R|_{\hat{\GC}_p}$ is homeomorphic to $S^2$ and the restriction 
	$\hat{\Gf}:\hat{\cX}_\R|_{\hat{\GC}_p}\to Z_\R$ is a homeomorphism.  
    In particular, $\{\del D_x\}_{x\in\GC_p}$ gives a foliation on 
	$Z_\R\setminus \{\text{\rm 2 points}\}$.  
  \item For each $p\in Z_\R$ and non zero $v\in T_pZ_\R$, 
    $\hat{\cX}_\R|_{\hat{\GC}_{p,v}}$ is homeomorphic to $S^2$ and the restriction 
	$\hat{\Gf}: \hat{\cX}_\R|_{\hat{\GC}_{p,v}} \to Z_\R$ is surjective. 
	Moreover, this is one-to-one distant from the curve $\hat{\Gf}^{-1}(p)$, hence 
    $\{(\del D_x\setminus\{p\})\}_{x\in\GC_{p,v}}$ gives a foliation on 
	$Z_\R\setminus \{p\}$.  
 \end{enumerate}
\end{Prop}

\begin{Rem}
 For distinguished points $p,q\in Z_\R\simeq\CP^1$, there are two families 
 of circles called ``Apollonian circles". 
 One of them is the family of the circles passing through $p,q$, which corresponds 
 to the space-like geodesic $\GC_{p,q}$. 
 The other family gives a foliation on $\CP^1\setminus\{p,q\}$, which corresponds 
 to a time-like geodesic and the foliation coincides with the one given in 
 {\it 1} of Proposition \ref{prop:foliation_on_Z_R}. 
 The case {\it 2} 
 of Proposition \ref{prop:foliation_on_Z_R} corresponds to the degenerated case. 
\end{Rem}

\section{Perturbation of holomorphic disks}
\label{Section:Holomorphic_disks}

We now investigate in the deformation of holomorphic disks. 
For a complex manifold $A$ and its submanifold $B$, 
we call simply {\it holomorphic disk on} $(A,B)$ 
for a continuous map $(\bD,\del \bD)\rightarrow(A,B)$ 
which is holomorphic on the interior of $\bD=\{ z\in\C \, |\, |z|\leq 1\}$. 

As in the previous Section, we put $Z=\CP^1\times\CP^1$ 
and $Z_\R=\{(\y,\bar{\y}^{-1}) \, |\, \y\in\CP^1 \}$. 
We have the family of holomorphic disks $\{D_{(\l,t)}\}$ 
defined from the double fibration  (\ref{eq:standard_double_fibration}), 
and we call each $D_{(\l,t)}$ {\it the standard disk}. 
In this Section, we treat a small perturbation $N$ of $Z_\R$, and prove that 
there is a natural $(S^2\times \R)$-family of holomorphic disks on $(Z,N)$ 
each of which is close to some standard disk. 
From the general theory by LeBrun \cite{bib:LeBrun}, 
one can show that there exists real three-parameter family of 
holomorphic disks on $(Z,N)$ near each standard disk. 
We, however, use the method given in \cite{bib:LM02}
so that we can treat more detail. 

\vspace{2mm}
First of all, we recall the $C^k$-topology of the space of 
deformations of $Z_\R$ in $Z$. 
A small perturbation $N$ of $Z_\R$ can be written 
$$ N=\left\{ \left. \(\y,\overline{\p(\y)}^{-1}\) \, \right|\, \y\in\CP^1 \right\} $$
using an automorphism $\p : \CP^1\to\CP^1$ which is sufficiently close 
to the identity map. 
Let $\{A_i\}$ be finitely many compact subsets and $\{B_i\}$ be 
open subsets on $\CP^1$ with complex coordinates $\y_i$, 
which satisfy 
$A_i\subset B_i$, $\p(A_i)\subset B_i$ and $\cup_i A_i=\CP^1$. 
Then $\p$ is identified with 
a combination of functions $ (h_i)_i$ where $h_i\in C^k(A_i,\C)$ 
are defined by $\p(\y_i)=\y_i+h_i(\y_i)$. 
The $C^k$-topology of the set of deformations of $Z_\R$ in $Z$ is 
defined by the norm 
$$ \| \p \|_{C^k}= \sup_i \| h_i \|_{C^k(A_i)} $$
where $ \| h_i \|_{C^k(A_i)} $ is suprema on $A_i$ of absolute values of all 
partial derivatives of $h_i$ whose order is less than or equal to $k$. 
In particular, for a compact subset $A\in\CP^1$ 
which is contained in a coordinated open $B$ 
and which satisfies $\p(A)\subset B$, 
$\| h\|_{C^k(A)}$ is sufficiently small 
if $\p$ is sufficiently close to the identity 
where $\p(\y)=\y+h(\y)$. 

\begin{Lem} \label{lem:deformation_of_holo_disk}
 Fix a standard holomorphic disk $D=D_{(\l,t)}$. 
 If $N\subset Z$ is the image of any embedding $\CP^1\hookrightarrow Z$ 
 which is sufficiently close to the standard one in the $C^{k+l}$-topology 
 with $k,l\geq 1$, 
 then there is a real three-parameter $C^l$-family of holomorphic disks on $(Z,N)$ 
 each of which is  $L^2_k$ close to $D$. 
\end{Lem}
\begin{proof}
 Since there is a transitive action of $\PSL(2,\C)$ on the standard disks, 
 we can assume $(\l,t)=(0,0)$, i.e. 
 $$D=\left\{ (z,z)\in Z \, | \, z\in \bD \right\} $$
 where $\bD=\{ z\in\C \, | \, |z|\leq 1 \}. $ 
 If we put $A=\{ \y\in\C \, | \, \frac{1}{2}\leq |\y| \leq 2 \}$, then 
  $N$ can be written 
  $$\left\{ \left. \(\y \, , \, \overline{\(\y+h(\y)\)}^{-1}\) \in Z \, 
    \right| \, \y\in A \right\} $$ 
 near $\del D$ using $h \in C^{k+l}(A)$ 
 whose $C^{k+l}$-norm is sufficiently small. 

 Then a small perturbation of $\del D$ is given as the image of
 $$S^1\to N : \, \c \mapsto \( e^{i(\c+u(\c))}, 
  \left[ e^{-i(\c+\bar{u}(\c))} + \bar{h}\(e^{i(\c+u(\c))}\) \right]^{-1} \) $$
 where $u$ is a $\C$-valued function on $S^1=\R/2\pi\Z$. 
 Here we denoted $\bar{u}(\c)=\overline{u(\c)}$ and 
 $\bar{h}(\y)=\overline{h(\y)}$. 
 Then we define the maps 
 $\GF_i:L^2_k(S^1,\C)\times C^{k+l}(A,\C)\to L^2_k(S^1,\C)$  by 
 \begin{equation} \label{eq:GF_i} 
  \begin{aligned} {}
   [\GF_1(u,h)](\c) &=e^{i(\c+u(\c))}, \\ 
   [\GF_2(u,h)](\c) &=\left[ e^{-i(\c+\bar{u}(\c))} 
    + \bar{h}\(e^{i(\c+u(\c))}\) \right]^{-1}. 
	\end{aligned}
 \end{equation}
 For given $h$, we want to arrange $u\in L^2_k(S^1,\C)$ 
 so that $[\GF_i(u,h)](\c)$ extends holomorphically to 
 $\{|z|<0\}$ for $z=e^{i\c}$. 
 Taking the derivation $\GF_i$, we obtain
 \begin{equation} \label{eq:GF_i*}
  \begin{aligned} {}
  [\GF_{1 *(0,0)}(\dot{u},\dot{h})](\c) &= ie^{i\c} \dot{u}(\c), \\ 
  [\GF_{2 *(0,0)}(\dot{u},\dot{h})](\c) 
   &= i e^{i\c} \bar{\dot{u}}(\c) - e^{2i\c} \bar{\dot{h}}(e^{i\c}). 
  \end{aligned}
 \end{equation}

 Now, we introduce some bounded operators (c.f.\cite{bib:LM02}). Set 
 $$ \begin{aligned}
 L^2\!\!\downarrow \, &= \left\{ \sum_{l< 0} a_l e^{il\c} \ \left| \ 
   a_l\in\C, \ \sum_{l<0} | a_l |^2<\infty \right.  \right\}, \\
 L^2_k\!\!\downarrow \, &= \left\{ \sum_{l< 0} a_l e^{il\c} \ \left| \ 
   a_l\in\C, \ \sum_{l<0} l^{2k}| a_l |^2<\infty \right.  \right\}
   =L^2_k(S^1,\C) \cap L^2\!\!\downarrow,  
 \end{aligned} $$
 and define $\Pi: L^2_k(S^1,\C)\rightarrow L^2_k\!\!\downarrow$ by 
 $$ \Pi\( \sum_{l=-\infty}^\infty a_l e^{il\c} \) = 
   \sum_{l<0} a_l e^{il\c}. $$
 Similarly let us define $\pye : L^2_k(S^1,\C) \rightarrow \C$ by 
 $$ \pye\( \sum_{l=-\infty}^\infty a_l e^{il\c} \) = a_0. $$ 
 Then, for $k, l\geq 1$, we define a $C^l$-map 
 $$ \GF : L^2_k(S^1,\C) \times C^{k+l}(A,\C) \longrightarrow 
   L^2_k\!\!\downarrow \times L^2_k\!\!\downarrow \times 
   C^{k+l}(A,\C) \times \C\times\C\times\C $$
 $$ \GF=(\Pi\action\GF_1)\times(\Pi\action\GF_2)\times
    \Lye \times (\pye\action\GF_1)\times (\pye\action\GF_2) \times \xye, $$
 where
 $$ \Lye :  L^2_k(S^1,\C) \times C^{k+l}(A,\C) \longrightarrow C^{k+l}(A,\C) $$
 is the factor projection, and
 $$ \xye :  L^2_k(S^1,\C) \times C^{k+l}(A,\C) \longrightarrow \C $$ 
 is given by 
 $$ \xye(u,h) = \frac{1}{2\pi} \int_0^{2\pi} u(\c)d\c, $$
 i.e. $\xye(u,h)=\pye(u)$. 
 $\GF$ is $C^l$ since $\Pi$, $\Lye$, $\pye$ and $\xye$ are 
 all bounded linear operators, and its derivative is given by 
 $$ \GF_*=(\Pi\action\GF_{1*})\times(\Pi\action\GF_{2*})\times
   \Lye\times(\pye\action\GF_{1*})\times(\pye\action\GF_{2*})\times\xye. $$ 

 In particular, if we write $ \dot{u}(\c)=\sum_n u_n e^{in\c}$, 
 then we obtain 
 $$ \GF_{*(0,0)}\left[ \begin{array}{c} \dot{u} \\ \dot{h} \end{array} \right]
  = \left[ \begin{array}{c} 
     \sum_{n<0} i u_{n-1} e^{in\c} \\[1mm]
     i \sum_{n<0} \bar{u}_{1-n} e^{in\c} 
	   - \Pi\(e^{2i\c}\bar{\dot{h}}(e^{i\c}) \)\\[1mm]
	 \dot{h} \\ iu_{-1} \\ i\bar{u}_1- \pye\(e^{2i\c}\bar{\dot{h}}(e^{i\c})\) \\
	 u_0
   \end{array} \right]. $$
 Since $\GF_{*(0,0)}$ has bounded inverse, the Banach-space inverse 
 function theorem tells us that there is an open neighborhood 
 ${\mathfrak U}$ of $(0,0)\in L^2_k(S^1,\C)\times C^{k+l}(A)$ and open 
 neighborhood ${\mathfrak V}$ of $\text{{\bf 0}}\in 
 L^2_k\!\!\downarrow \times L^2_k\!\!\downarrow \times 
 C^{k+l}(A,\C) \times \C\times\C\times\C$ such that 
 $\GF|_{\mathfrak U} : {\mathfrak U}\rightarrow {\mathfrak V}$ 
 is a $C^l$-diffeomorphism. 

 Hence, for given $h$, we obtain complex three-parameter $C^l$-family 
 of holomorphic disks defined from 
 $(u,h)= \GF^{-1}(0,0,h,\a_1,\a_2,\b)$, where $\a_1,\a_2,\b$ are 
 small complex numbers. 
 It contains, however, real three-dimensional ambiguity which comes from the 
 disk automorphism. 
 To kill this ambiguity, it is enough to use the inverse of 
 \begin{equation} \label{eq:three_para_families}
 (0,0,h,\a,-\a,i\b)\, \in \, L^2_k\!\!\downarrow \times L^2_k\!\!\downarrow 
 \times C^{k+l}(A,\C) \times \C\times\C\times\C, 
 \end{equation}
 by $\GF$ for $(\a,\b)\in\C\times\R$ which is sufficiently close to $(0,0)$. 
 Now the statement follows since $\| h\|_{C^{k+l}(A)}$ is sufficiently small if 
 $N$ is sufficiently close to $Z_\R$. 
\end{proof}

\begin{Rem} \label{rem:local_family}
 1. \, 
 Let $\GD$ be any holomorphic disk on $(Z,N)$ constructed as above lemma. 
 Then $\GD$ intersects with $N$ only on the boundary $\del \GD$. 
 Actually, let $\bD\to Z: z\mapsto (\vp_1(z),\vp_2(z))$ be the map 
 corresponding to $\GD$ and denote 
 $N=\{ (\y\, ,\, \overline{\p(\y)}^{-1}) 
   \, | \, \y\in\CP^1 \}.$ 
 Notice that $\y\mapsto\overline{\p(\y)}^{-1}$ maps $\vp_1(\del \bD)$ to 
 $\vp_2(\del \bD)$ and maps the interior of $\vp_1(\bD)$ to the outside of 
 $\vp_2(\bD)$. 
 Suppose that there is an interior point $z\in \bD$ such that 
 $\vp_2(z)=\overline{\p(\vp_1(z))}^{-1}$. Then 
 $\vp_1(z)$ is contained in the interior of $\vp_1(\bD)$, and 
 $\overline{\p(\vp_1(z))}^{-1}$ is out side of $\vp_2(\bD)$. 
 This is a contradiction. 

 2. \,  
 We can take $\GV$ so that 
\begin{equation} \label{eq:V_as_a_product}
 \begin{aligned}
 \GV&=\GV_1\times\GV_2\times\GW\times V_1 \times V_2 \times V_3, \\
 \GW&=\left\{ \left. h\in C^{k+l}(A,\C) \, \right| \, 
      \| h \|_{C^{k+l}(A)} < \e_0 \right\}, 
	 \end{aligned}
 \end{equation}
 where $\GV_i\subset L^2_k\!\!\downarrow$ and 
 $V_i\subset \C$ are small open sets and $\e_0>0$ is a constant. 
 This notation is used in the following arguments. 
\end{Rem}

\vspace{2mm}
Next we want to prove that, if $N$ is sufficiently close to $Z_\R$, 
then the method of Lemma \ref{lem:deformation_of_holo_disk} 
works for {\it all} standard disks at once. 
Then we need uniform estimate of the deformation $N$ of $Z_\R$ 
among all standard disks. 
In the LeBrun-Mason's case \cite{bib:LM02,bib:LM05}, 
the parameter spaces of holomorphic disks are compact and homogeneous, 
so the uniform estimate is automatically deduced from the local estimate. 
In our case, however, the parameter space is non-compact space $S^2\times\R$, 
so we need more detail arguments. 
For this purpose, it is enough to show that
the deformations of the disks are ``tame" as in the following lemma
on the neighborhood of the boundary of the parameter space. 

\begin{Lem} \label{lem:deformation_near_bdry}
 Let $\{D_{(\l,t)}\}$ be the standard disks. 
 Suppose $N\subset Z$ is sufficiently close to $Z_\R$ in the $C^{k+l}$-topology. 
 Then the three-parameter family of holomorphic disks on 
 $(Z,N)$ near $D_{(\l,t)}$ always exists 
 for each $(\l,t)\in\CP^1\times\R$ with $t\gg 0$.  
\end{Lem}

\begin{proof}
 It is enough to consider the case $\l=0$. 
 We fix a small constant $c>0$ and 
 let $B_c=\{z\in\C \, | \, |z|< c\}$. 
 Notice that the compact subset $B_{c}\times \CP^1\subset Z$ 
 contains all holomorphic disks of the form $D_{(0,t)}$ if $e^t >2c^{-1}$. 
 We can write 
 \begin{equation} \label{eq:N_in_Bc_times_CP^1}
 N\cap (B_c\times \CP^1) = 
 \left\{ \left. \(\y, \overline{(\y+h(\y))}^{-1}\) \, \right| \, \y\in B_c \right\}
 \end{equation}
 using $h\in C^{k+l}(B_c,\C)$. 
 We claim that if $\|h\|_{C^{k+l}(B_c)}<\frac{\e_0}{4\sqrt{2}}$,  
 then the three-parameter family of 
 holomorphic disks on $(Z,N)$ near $D_{(0,t)}$ exists for all $e^t >2c^{-1}$. 
 Here $\e_0$ is the constant defined in (\ref{eq:V_as_a_product}). 

 Now we show that it is enough to prove the case when $h(0)=0$ and 
 $\|h\|_{C^{k+l}(B_c)}<\frac{\e_0}{2\sqrt{2}}$. 
 Actually, if we change the coordinate $(\y_1,\bar{\y}_2^{-1})\in Z$ to  
 $(\x_1,\bar{\x}_2^{-1})$ by the relation 
 $$\x_1=\y_1,\quad \x_2=\y_2+h(0), $$ 
 then we can write 
 $$ N\cap (B_c\times \CP^1) = 
 \left\{ \left. \(\x, \overline{(\x+g(\x))}^{-1}\) \, \right| \, \x\in B_c \right\}$$
 using $g(\x)=h(\x)-h(0)$. 
 Here we obtain $\|g\|_{C^{k+l}(B_c)}<\frac{\e_0}{2\sqrt{2}}$ since 
 $$ \begin{aligned}
  \sup_{\x\in B_c}|g(\x)| & < \sup_{\x\in B_c}|h(\x)|+ |h(0)|
    <\frac{\e_0}{2\sqrt{2}}, \\
  \sup_{\x\in B_c}|Dg(\x)| & = \sup_{\x\in B_c}|D h(\x)| < \frac{\e_0}{4\sqrt{2}}, 
  \end{aligned} $$ 
 where $D$ is any partial derivative of the degree less than or equal to $k+l$. 
 Hence, if we rewrite $h$ instead of $g$, 
 we can assume $h(0)=0$ and $\|h\|_{C^{k+l}(B_c)}<\frac{\e_0}{2\sqrt{2}}$

 We denote $r=e^t$ from now on. 
 A small perturbation of $\del D_{(0,t)}$ is given as the image of 
 $$S^1\to N : \, \c \mapsto \( r^{-1}e^{i(\c+u(\c))}, 
  \left[ r^{-1}e^{-i(\c+\bar{u}(\c))} + \bar{h}\(r^{-1}e^{i(\c+u(\c))}\) 
  \right]^{-1} \), $$
 where $u$ is a $\C$-valued function on $S^1$. 
 
 Let $A^r=\{z\in\C \, | \, \frac{r^{-1}}{2} \leq |z| \leq 2r^{-1} \}$ and $A=A^1$, 
 then $A^r$ is a compact subset of $B_c$ if $r>2c^{-1}$. 
 We define the maps 
 $\GF^r_i:L^2_k(S^1,\C)\times C^{k+l}(A^r,\C)\to L^2_k(S^1,\C)$  by 
 $$ \begin{aligned} {}
   [\GF^r_1(u,h)](\c) &=r^{-1}e^{i(\c+u(\c))}, \\ 
   [\GF^r_2(u,h)](\c) &=\left[ r^{-1}e^{-i(\c+\bar{u}(\c))} 
    + \bar{h}\(r^{-1}e^{i(\c+u(\c))}\) \right]^{-1}. 
	\end{aligned} $$
 Putting $h^r(z)=r\, h(r^{-1}z)$, we obtain 
 \begin{equation} \label{eq:GF^e_i}
  [\GF^r_1(u,h)](\c)=r^{-1}[\GF_1(u,h^r)](\c), \qquad  
  [\GF^r_1(u,h)](\c)=r [\GF_1(u,h^r)](\c), 
 \end{equation}
 where $\GF_i$ is the map given by (\ref{eq:GF_i}). 
 Notice that the map $\r^r:h\mapsto h^r$ gives 
 an isomorphism of Banach spaces $ C^{k+l}(A^r,\C)\to C^{k+l}(A,\C) $. 

 Similar to the proof of lemma \ref{lem:deformation_of_holo_disk}, we define 
 $$ \GF^r : L^2_k(S^1,\C) \times C^{k+l}(A^r,\C) \longrightarrow 
   L^2_k\!\!\downarrow \times L^2_k\!\!\downarrow \times 
   C^{k+l}(A^r,\C) \times \C\times\C\times\C $$
 $$ \GF^r=(\Pi\action\GF^r_1)\times(\Pi\action\GF^r_2)\times
    \Lye \times (\pye\action\GF^r_1)\times (\pye\action\GF^r_2) \times \xye, $$
 where $\Lye$ is the projection. 
 Then we can relate $\GF^r$ with $\GF$ in the following way. 
 Let $m(r)$ be multiplication of $r$ on $L^2_k\!\!\downarrow$ or $\C$,  
 then we obtain the following commutative diagram 
 \begin{equation} \label{eq:comm_diag_of_expand_and_shrink}
 \xymatrix{  L^2_k(S^1,\C) \times C^{k+l}(A^r,\C) \ar[rr]^{\GF^r\qquad\quad} 
  \ar[d]^{\text{id}\times \r^r} & &
   L^2_k\!\!\downarrow \times L^2_k\!\!\downarrow \times 
   C^{k+l}(A^r,\C) \times \C\times\C\times\C \ar[d]^{\Phi^r} \\
   L^2_k(S^1,\C) \times C^{k+l}(A,\C) \ar[rr]^{\GF\qquad\quad}  & &
   L^2_k\!\!\downarrow \times L^2_k\!\!\downarrow \times 
   C^{k+l}(A,\C) \times \C\times\C\times\C } 
 \end{equation}
 where $ \Phi^r=m(r)\times m(r^{-1})\times\r^r\times 
   m(r)\times m(r^{-1})\times\text{id}. $
 Notice that the vertical arrows in the above diagram are isomorphisms, 
 and that the restriction $\GF|_\GU:\GU\to\GV$ gives $C^l$-diffeomorphism 
 from the proof of lemma \ref{lem:deformation_of_holo_disk}. 
 Hence the restriction 
 $$\GF^r : (\text{id}\times\r^r)^{-1}(\GU)\longrightarrow (\Phi^r)^{-1}(\GV) $$
 is $C^l$-diffeomorphism. 
 If we take $\GV$ to be the product as in (\ref{eq:V_as_a_product}), then 
 $$ (\Phi^r)^{-1}(\GV)= r^{-1}\GV_1 \times r\GV_2 
  \times (\r^r)^{-1}(\GW) \times r^{-1}V_1 \times rV_2 \times V_3. $$

 We want to show that $h|_{A^r}\in(\r^r)^{-1}(\GW)$, 
 or equivalently $\|h^r\|_{C^{k+l}(A)}<\e_0$, for all $r>2c^{-1}$.  
 Let $x,y$ be the real coordinate such that $\y=x+iy$, and let 
 $D=\del^m/\del x^j \del y^{m-j}$ be a derivation of degree $m\leq l+k$, 
 then we obtain 
 $$ Dh^r(\y)= r^{1-m} D h(r^{-1}\y).$$ 
 Hence 
 $$ \sup_{\y\in A}|Dh^r(\y)| \leq r^{1-m}\sup_{\y\in A}|Dh(r^{-1}\y)|
 \leq r^{1-m}\sup_{\z\in A^r}|Dh(\z)| <\frac{\e_0}{2\sqrt{2}}r^{1-m}<\e_0, $$  
 if $m\geq 1$. For $m=0$, notice that 
 $$ \begin{aligned}
  |h(\y)| & \leq \int_0^1 \left| \frac{dh}{dt}(t\y) \right| dt
    \leq \int_0^1 \left| \frac{\del h}{\del x}(t\y) \right| |x| dt
          +\int_0^1 \left| \frac{\del h}{\del y}(t\y) \right| |y| dt \\
	& < \frac{\e_0}{2\sqrt{2}}(|x|+|y|) < \frac{\e_0}{2} |\y|, 
  \end{aligned} $$ 
 hence we obtain 
 $$ \sup_{\y\in A}|h^r(\y)| = r \sup_{\y\in A}|h(r^{-1}\y)|
    = r \sup_{\z\in A^r}|Dh(\z)| < \frac{r\e_0}{2} \sup_{\z\in A^r}|\z| = \e_0. $$
 In this way, we have obtained $\|h^r(\y)\|_{C^{k+l}(A)}<\e_0$ for all $r>2c^{-1}$. 
\end{proof}

\begin{Rem} \label{Rem:another_limit_of_holo_disk}
 Lemma \ref{lem:deformation_near_bdry} also holds for  $t\ll 0$. 
 Exchange the role of factors of $Z=\CP^1\times \CP^1$ and change $t$ with $-t$ 
 to prove this case. 
\end{Rem}

From Lemma \ref{lem:deformation_of_holo_disk}
and \ref{lem:deformation_near_bdry}, we obtain the following statement. 

\begin{Prop} \label{Prop:perturbation_of_all_disks}
 If $N\subset Z$ is the image of any embedding $\CP^1\hookrightarrow Z$ 
 which is sufficiently close to the standard one in the $C^{k+l}$-topology with 
 $k,l\geq 1$, 
 then there is a $C^l$ family of holomorphic disks on $(Z,N)$ 
 each of which is $L^2_k$ close to some standard disk on $(Z,Z_\R)$. 
\end{Prop}

We will strengthen this statement in 
Proposition \ref{Prop:family_of_perturbed_holo_disks}. 

\section{The double fibration}
\label{Section:The_double_fibration}

In this section, 
we investigate in some properties for the family of holomorphic disks 
constructed in Section \ref{Section:Holomorphic_disks}. 
We continue to use the notations $\GF,\GF_i,\GU,\GV$ and so on. 

For each $h\in C^{k+l}(A,\C)$, 
we define $C^l$-maps $\Xi^h, F^h_i : U \to L^2_k(S^1,\C)$ so that 
$$ \begin{aligned}
   \(\Xi^h(\a,\b),h\) &=\GF^{-1}(0,0,h,\a,-\a,i\b), \\ 
   F^h_1(\a,\b)(e^{i\c}) &= \GF_1\(\Xi^h(\a,\b),h\)(\c) = 
     \exp i\left\{\c+\Xi^h(\a,\b)(\c)\right\}, \\ 
   F^h_2(\a,\b)(e^{i\c}) &= \GF_2\(\Xi^h(\a,\b),h\)(\c), 
	 \end{aligned} $$
where $U\subset \C\times\R$ is a small open neighborhood of $(0,0)$
depending on $h$. 
By definition, $F^h_i(\a,\b)(z)$ extend to holomorphic functions on 
$\bD=\{|z|\leq 1\}$, and satisfy $F^h_1(\a,\b)(0)=\a$
and $F^h_2(\a,\b)(0)=-\a$. 
If we expand 
\begin{equation} 
 \Xi^h(\a,\b)(\c)=\sum_k \Xi^h(\a,\b)_k e^{ik\c}, 
\end{equation}
then we obtain $\Xi^h(\a,\b)_0=i\b$ by definition. 
Notice that we can also define the derivations $\Xi^h_*$ and $F^h_{i\, *}$ 
which satisfy 
$$ \begin{aligned}
   \(\Xi^h_*(\dot\a,\dot\b),0\) &=\GF^{-1}_*(0,0,0,\dot\a,-\dot\a,i\dot\b), \\ 
   F^h_{1\, *}(\dot\a,\dot\b)(e^{i\c}) 
    &=\GF_{1\, *}\(\Xi^h(\dot\a,\dot\b),0\)(\c)  
       =i \, F^h_1(e^{i\c}) \, \Xi^h_*(\dot\a,\dot\b)(\c), \\ 
   F^h_{2\, *}(\dot\a,\dot\b)(e^{i\c}) 
    &= \GF_{2\, *}\(\dot\Xi^h(\dot\a,\dot\b),0\)(\c), 
	\end{aligned} $$
$$ F^h_{1\, *}(\dot\a,\dot\b)(0)=\dot\a, \qquad
  F^h_{2\, *}(\dot\a,\dot\b)(0)=-\dot\a, \quad \text{and} \quad 
  \Xi^h_*(\dot\a,\dot\b)_0=i\dot\b. $$

Let $N\subset Z$ be the image of any embedding $\CP^1\hookrightarrow Z$ 
which satisfies Proposition \ref{Prop:perturbation_of_all_disks}. 
Let us denote $\GB^N_{(\a,\b)}$ for the holomorphic disk on $(Z,N)$ 
which corresponds to the element $(0,0,h,\a,-\a,i\b)\in\GV$ 
in the notation of the proof of Lemma \ref{lem:deformation_of_holo_disk}. 
Then 
\begin{equation} \label{eq:GB} 
 \GB^N_{(\a,\b)}=\left\{ \left. \( F^h_1(\a,\b)(z),F^h_2(\a,\b)(z) \) \in Z 
   \, \right| \, z\in\bD \right\}, 
\end{equation}
and $\{\GB^N_{(\a,\b)}\}_{(\a,\b)\in U}$ gives the three-parameter family of 
holomorphic disks each of which is $L^2_k$-close to
the standard disk $D_{(0,0)}$. 
Notice that $\GB^N_{(\a,\b)}$ passes through $(\a,-\a)$ when $z=0$, 
hence, for fixed $\a$, $\{\GB^N_{(\a,\b)}\}_\b$ defines a one-parameter 
family of holomorphic disks which pass through $(\a,-\a)$. 

\vspace{2mm}
In the standard case, the following statement holds. 
\begin{Prop}
 $\GB^{\Z_\R}_{(\a,\b)}$ coincides with the standard disk $D_{(\a,\b)}$.
\end{Prop}
\begin{proof}
 Since the disk 
 $$ \GB^{\Z_\R}_{(\a,\b)}=\left\{ \left. \(F^0_1(\a,\b)(z), F^0_2(\a,\b)(z) \) \in Z 
   \, \right| \, z\in\bD \right\} $$ 
 coincides with one of the standard disks near $D_{(0,0)}$, 
 there is a unique element $(\l,t)\in\CP^1\times\R$ near $(0,0)$ such that 
 \begin{equation} \label{eq:interpretation_of_Mobius_transf}
   F^0_1(\a,\b)(e^{i\c})= \exp i\( \c+ \Xi^0(\a,\b)(\c)\) 
   = \frac{e^{i\c}+e^t\l}{-\bar{\l}e^{i\c}+e^t}. 
 \end{equation}
 Then we obtain $\a=F^0(\a,\b)(0)=\l$.  
 On the other hand, taking the derivation of (\ref{eq:interpretation_of_Mobius_transf}), 
 we obtain 
 $$ i\,\Xi^0_*(\dot\a,\dot\b)(\c) = 
   \frac{(\dot{\l}+\l\dot{t})e^{-(i\c-t)}}{1+\l e^{-(i\c-t)}} + 
   \frac{e^{i\c-t}\bar{\dot{\l}}-\dot{t}}{1-\bar{\l}e^{i\c-t}}. $$
 If we expand the right hand side and compare the constant terms, then we find 
 $$ i\dot{\b}=\Xi^0_*(\dot\a,\dot\b)_0 = i\dot{t}. $$
 On the other hand, it is easy to see that $t=\b$ when $\a=0$. 
 Hence $(\l,t)=(\a,\b)$ for each $(\a,\b)\in U$.
\end{proof}

Let $M$ be the parameter space 
of the family of holomorphic disks on $(Z,N)$ constructed 
in Proposition \ref{Prop:perturbation_of_all_disks}. 
Then $M$ has natural structure of real 3-manifold and 
we can take a coordinate system on $M$ in the following way. 
For each $(\l,t)\in\CP^1\times\R$, chose an element $T=T{(\l,t)}\in\PSL(2,\C)$ 
such that $T_*(D_{(\l,t)})=D_{(0,0)}$, where
$\{D_{(\l,t)}\}$ is the standard disks. 
Let $U^T\subset\C\times\R$ be an open neighborhood of $(0,0)$ 
such that $\GB^{T_*(N)}_{(\a,\b)}$ are defined. 
Then $\left\{T^{-1}_*\GB^{T_*(N)}_{(\a,\b)}\right\}_{(\a,\b)\in U^T}$ gives 
the family of holomorphic disks on $(Z,N)$ each of which is 
close to $D_{(\l,t)}$, 
and $\{U^{T(\l,t)}\}$ gives a coordinate system on $M$. 

Using above coordinates, we prove the following Lemma. 
\begin{Lem} \label{Lem:R_family_pass_through_q_in_Q}
 Suppose $N\subset Z$ is sufficiently close to $Z_\R$ so that 
 Proposition \ref{Prop:perturbation_of_all_disks} holds, 
 and consider the 
 constructed family of holomorphic disks on $(Z,N)$. 
 Then, for each $q=(\a,-\a)\in Z$, 
 there is an $\R$-family of holomorphic disks 
 each of which passes through $q$. 
 Moreover there is a natural compactification of this family and 
 the boundary points $\pm\infty$ correspond to marked $\CP^1$. 
\end{Lem}
\begin{proof}
 We can assume $\a=0$. 
 Take any $t$ so that $|t|$ is sufficiently small, 
 and consider the standard disk $D_{(0,t)}$. 
 If we define $T\in\PSL(2,\C)$ by 
 $$ T=\left[ \begin{array}{cc} e^{\frac{t}{2}} & \\ &e^{-\frac{t}{2}} 
  \end{array} \right], $$
 then $T_*(\y_1,\y_2)=(e^t\y_1,e^{-t}\y_2)$ and $T_*(D_{(0,t)})=D_{(0,0)}$. 
 Now $\left\{T^{-1}_*\GB^{T_*(N)}_{(0,\b')}\right\}_{\b'\in V} $ gives a 
 one-parameter family of holomorphic disks on $(Z,N)$ each of which is 
 close to $D_{(0,t)}$ and pass through $(0,0)$. 
 Here we denoted $V=\{ \b'\in\R \, | \, (0,\b')\in U^T \}$. 

 Since $|t|$ is small, there is an open set $V'\subset V$ such that
 $T^{-1}_*\GB^{T_*(N)}_{(0,\b')}$ is sufficiently close to $D_{(0,0)}$ 
 for all $\b'\in V'$. 
 Hence, for each $\b'\in V'$, there is a unique $(\a,\b)$ such that 
 \begin{equation} \label{eq:extension_of_R-family}
   T^{-1}_*\GB^{T_*(N)}_{(0,\b')} = \GB^N_{(\a,\b)}. 
 \end{equation}

 Now $N$ and $T_*(N)$ can be written locally 
  \begin{equation} \label{eq:N_and_T_*(N)}
  N: \left\{ \left. \(\y,\overline{\(\y+h(\y)\)}^{-1}\) \ \right| 
      \, \y \in A \right\}, \qquad 
   T_*(N): \left\{ \left. \(\y,\overline{\(\y+h^T(\y)\)}^{-1}\) \ \right| 
      \, \y \in A \right\} 
  \end{equation}
 using $C^{k+l}$-function $h$ which is defined on a neighborhood of 
 $A=\{z\in\C\, |\, \frac{1}{2}\leq |z| \leq 2 \}$. 
 Here we denoted $h^T=ThT^{-1}$. 
 Then (\ref{eq:extension_of_R-family}) is equivalent to 
 $$ e^{-t} F^{h^T}_1(0,\b')(z)=F^{h}_1(\a,\b)(z) \qquad \text{on} \quad z\in\bD. $$
 Evaluating $z=0$, we obtain $\a=0$. 
 Moreover, this is also equivalent to 
 $$ it+\Xi^{h^T}(0,\b')(\c) = \Xi^h(\a,\b)(\c) \qquad \text{on} \quad \c\in S^1. $$
 Comparing the constant terms for $e^{i\c}$, we obtain $\b=\b' +t$. 
 Hence (\ref{eq:extension_of_R-family}) is equivalent to $(\a,\b)=(0,\b'+t)$. 
 So the one-parameter family $\{\GB^N_{(0,\b)}\}_{(0,\b)\in U}$ 
 extends to 
 $$ \left\{ \b\in\R \, \left| \ (0,\b)\in U \, \text{or} \, (0,\b-t)\in U^T 
   \right. \right\}$$ 
 by putting $\GB^N_{(0,\b)}= T^{-1}_* \GB^{T_*(N)}_{(0,\b-t)}$. 
 In this way, we can define the one-parameter family $\{\GB^N_{(0,\b)}\}$ 
 for all $\b\in\R$. 
 
 The statement of the compactification is obtained 
 from Lemma \ref{lem:deformation_near_bdry} and its proof. 
 Indeed, in the notation of (\ref{eq:N_in_Bc_times_CP^1}), 
 if we take the limit $t\to\infty$, then we obtain the family of holomorphic disks 
 on $(Z,N)$ whose limit degenerates to 
 $\{0\} \times \CP^1$ marked at $(0,\overline{h(0)}^{-1})$. 
 As we explained at Remark \ref{Rem:another_limit_of_holo_disk}, 
 we also obtain another marked $\CP^1$ by $t\to -\infty$. 
\end{proof}

Now the following statement is easily proved. 
\begin{Prop} \label{Prop:family_of_perturbed_holo_disks}
 If $N\subset Z$ is the image of any embedding $\CP^1\hookrightarrow Z$ 
 which is sufficiently close to the standard one in the $C^{k+l}$-topology with 
 $k,l\geq 1$, 
 then there is a $C^l$ family of holomorphic disks on $(Z,N)$ 
 parametrized by $S^2\times\R$
 which satisfies the following properties: 
 \begin{itemize}
  \item each disk is $L^2_k$-close to some standard disk, 
  \item there is a natural compactification of the family such that 
   the compactified family is parameterized by $S^2\times I$, 
   and each boundary point on $S^2\times I$ 
   corresponds to a marked $\CP^1$ embedded in $(Z,N)$, 
 \end{itemize}
 where $I=[-\infty,\infty]$ is the compactification of $\R$. 
\end{Prop}
\begin{proof}
 Let $Q=\{(\l,-\l)\in Z\, |\, \l\in\CP^1\}$. For each $q\in Q$, there is an 
 $\R$-family of holomorphic disks constructed in 
 Lemma \ref{Lem:R_family_pass_through_q_in_Q}. 
 Since these families varies continuously, we obtain the family of holomorphic disks 
 parametrized by $Q\times\R\simeq S^2\times\R$. 
 The statement for the compactification is obvious from 
 Lemma \ref{Lem:R_family_pass_through_q_in_Q}. 
\end{proof}

For each $(\l,t)\in\CP^1\times\R$, we define 
$$ \GD_{(\l,t)}=T^{-1}_*\GB^{T_*(N)}_{(0,0)} $$ 
where $T=T(\l,t)\in \PSL(2,\C)$ is 
an element which satisfies $T_*(D_{(\l,t)})=D_{(0,0)}$. 
Then we obtain the continuous map 
$j:\CP^1\times\R\to M : (\l,t) \mapsto \GD_{(\l,t)}$. 
Moreover, we can prove that $j$ is an isomorphism in the following way.
For each constructed holomorphic disk $\GD$ on $(Z,N)$, 
we can choose $(\l,t)$ and $T=T(\l,t)$ so that 
$\GD=T^{-1}_*\GB^{T_*(N)}_{(0,\b)}$. 
Here $\l$ is uniquely defined so that the center of $\GD$ is $(\l,-\l)$. 
Then $\GD=\GD_{(0,\b+t)}$ from Lemma \ref{Lem:R_family_pass_through_q_in_Q} 
and its proof, so $j$ is surjective. 
The injectivity and the continuity of $j^{-1}$ is also deduced from above procedure 
of choosing $(\l,t)$, hence $j$ is isomorphism. 

\vspace{2mm}
Let us construct the double fibration. 
Let $U\subset\CP^1\times\R$ be a sufficiently small neighborhood of $(0,0)$. 
For each $(\l,t)\in U$, we define $T=T(\l,t)\in\PSL(2,\C)$ by 
 $$ T=\frac{1}{e^{-\frac{t}{2}}\sqrt{1+|\l|^2}} \left[ 
  \begin{array}{cc}  1&-e^t\l \\ \bar{\l}&e^t \end{array} \right], $$
then we obtain $T_*(D_{(\l,t)})=D_{(0,0)}$. 
Introducing $C^{k+l}$-function $h$ and $h^T$ 
similarly to (\ref{eq:N_and_T_*(N)}), 
we define a map $\Gf:U\times\bD\to Z$ by 
$$ \Gf(\l,t;z)= T_*^{-1}\(F^{h^T}_1(z), F^{h^T}_2(z) \). $$
Then $\Gf$ is $C^l$ for $(\l,t)$ and $C^{k-1}$ for $z$, and we obtain 
$\GD_{(\l,t)}=\{\Gf(\l,t;z)\in Z \, | \, z\in\bD\}.$ 

Constructing similar map for each neighborhood of $\CP^1\times\R$, 
and patching them, we obtain the double fibration 
\begin{equation} \label{eq:perturbed_double_fibration}
 \xymatrix{
 & (\cX_+, \cX_\R) \ar[dl]_{\varpi} \ar[dr]^{\Gf} & \\
 M\simeq\CP^1\times\R && (Z,N)}
\end{equation}
where $\varpi$ is a disk bundle. 
By the construction, $\cX_+$ is the same disk bundle as the standard case. 
In particular, we obtain that $c_1(\cX_\R)=2$ along each fiber of $\varpi$ 
and that $\varpi$ is $C^\infty$. 

\begin{Lem} \label{Lem:non_vanishing_of_tangent_vec_of_circles}
 Let $N\subset Z$ be the image of any embedding $\CP^1\hookrightarrow Z$ 
 which is sufficiently close to the standard one in the $C^{k+l}$-topology with 
 $k, l \geq 1$, and consider the double fibration (\ref{eq:perturbed_double_fibration}). 
 Then $f_*({\bf v})\neq 0$ for each non zero vector $\text{\bf v}\in T\cX_\R$ 
 such that $\varpi_*({\bf v})=0$. 
\end{Lem}
\begin{proof}
 For each $(u,h)\in L^2_k(S^1,\C)\times C^{k+l}(A,\C)$, we have 
 $$\frac{d}{d\c}[\GF_1(u,h)](\c)=\frac{d}{d\c}e^{i(\c+u(\c))}
    =e^{i(\c+u(\c))}(i+iu'(\c)), $$
 so this dose not vanish if $\|u\|_{L^2_1}$ is sufficiently small. 
 Hence, by shrinking $\GU$ and $\GV$ smaller if we need, the statement holds 
 for $\text{\bf v}\in \ker\varpi_*$ over $U\subset M$ 
 where $U$ is the neighborhood introduced above. 

 Now, recall the diagram (\ref{eq:comm_diag_of_expand_and_shrink}) in the 
 proof of Lemma \ref{lem:deformation_near_bdry}. 
 Notice that the $L^2_k(S^1,\C)$ component dose not change by the vertical arrow, 
 so we can estimate $u\in L^2_k(S^1,\C)$ uniformly so that 
 $\frac{d}{d\c}[\GF_1^r(u,h)](\c)$ does not vanish for all $r$. 
 Hence the statement holds for all $\text{\bf v}\in \ker\varpi_*$. 
\end{proof}

By Lemma \ref{Lem:non_vanishing_of_tangent_vec_of_circles}, 
we can define the lift $\tilde\Gf$ of $\Gf$ by 
$\tilde\Gf: \cX_R\to S(TN) : u\mapsto [\Gf_*(\text{\bf v}_u)]$. 
Here {\bf v} is a nowhere vanishing vertical vector field, i.e. $\varpi_*({\bf v})=0$,  
whose orientation matches to 
the complex orientation of the fiber of $\varpi:\cX_+\to M$. 
The next Proposition is the perturbed version of 
Proposition \ref{Prop:properties_over_S_t}. 
\begin{Prop} \label{Prop:properties_over_S_t_general_version}
 Let $N\subset Z$ be the image of any embedding $\CP^1\hookrightarrow Z$ 
 which is sufficiently close to the standard one in the $C^{k+l}$-topology with 
 $l\geq 1, k\geq 2$. 
 Consider the double fibration (\ref{eq:perturbed_double_fibration}), 
 let $S_t=\CP^1\times \{t\} \subset M$, 
 and  let $\Gf_t$ and $\tilde{\Gf}_t$ be the restriction of $\Gf$ and $\tilde{\Gf}$ on 
 $\varpi^{-1}(S_t)$ respectively. Then, for each $t\in \R$, 
 \begin{enumerate}
  \item $\Gf_t:(\cX_+\setminus\cX_\R)|_{S_t}\to Z\setminus N$  
    is diffeomorphic, 
 \item $\tilde{\Gf}_t:\cX_\R|_{S_t}\to S(TN)$ is diffeomorphic, 
 \item $\Gf_t:\cX_\R|_{S_t}\to N$ is an $S^1$-fibration such that 
    each fiber is transverse to the vertical distribution of $\varpi: \cX_\R\to M$, 
 \end{enumerate}
 In particular, $\{\GD_{(\l,t)}\}_{\l\in\CP^1}$ 
 gives a foliation on $Z\setminus N$ for each $t\in\R$. 
\end{Prop}
\begin{Rem}
 From {\it 2} above, the following holds:  
 for each $t\in\R$, $p\in N$ and non zero $v\in T_p N$, there is a unique 
 $x\in S_t$ such that $p\in\del\GD_x$ and $v\parallel \GD_x)$. 
\end{Rem}

\begin{proof}[Proof of \ref{Prop:properties_over_S_t_general_version}]
 Since $S_t$ is compact and $\Gf$ is $C^1$-close to the standard case, 
 we can assume the derivation of $\Gf_t$ to be non zero everywhere 
 by shrinking $\GW$ smaller if we need. 
 Here $\GW$ is the open set defined in Remark \ref{rem:local_family}. 
 Notice that we can arrange so that this property holds for all $t\in\R$ at once 
 by Lemma \ref{lem:deformation_near_bdry} and its proof. 
 Thus $\Gf_t$ gives proper local diffeomorphism on $(\cX_+\setminus\cX_\R)|_{S_t}$, 
 and this is actually diffeomorphism since $\Gf_t$ is close to the standard case. 
 
 By the similar argument for the lift $\tilde{\Gf} : \cX_\R|_{S_t}\to S(TN)$, 
 we obtain the property {\it 2}. 
 If there are $x\in S_t$ and $p\in N$ such that 
 $\varpi^{-1}(x)$ and $\Gf_t^{-1}(p)$ are not transversal at $u\in\cX_\R$, 
 then $(\Gf_t)_*(\text{\bf v}_u)=0$. 
 This contradicts to Lemma \ref{Lem:non_vanishing_of_tangent_vec_of_circles}, 
 Hence {\it 3} holds. 
\end{proof}

From Proposition \ref{Prop:family_of_perturbed_holo_disks}, 
we obtain the natural compactification of $\varpi$ and $\Gf$ 
which gives the following double fibration: 
\begin{equation} \label{eq:compactified_double_fibration_perturbed_case} 
 \xymatrix{
  & (\hat{\cX}_+, \hat{\cX}_\R) \ar[dl]_{\hat\varpi} \ar[dr]^{\hat\Gf} & \\
  \hat{M} && (Z,N)} 
\end{equation}
which is studied in Section \ref{Section:Construction_of_E-W_spaces}. 

\vspace{2mm}
For the last of this section, we prove the following technical lemma 
which enables us to prove the non-degeneracy of the induced conformal structure. 
Let us denote $\GC_p=\varpi\action\Gf^{-1}(p)=\{x\in M\, |\, p\in\GD_x \}$ 
for each $p\in Z\setminus N$, then $\GC_p$ is an embedded $\R$ in $M$ from 
Proposition \ref{Prop:properties_over_S_t_general_version}. 
Notice that $\GC_p$ is a closed subset in $M$ since it 
connects two boundaries of $\hat{M}$. 

\begin{Lem} \label{lem:non-degeneracy_of_time-like_vectors}
 Let $x\in M$, 
 then there are two points $p_1,p_2\in \GD_x\setminus\del\GD_x$ 
 such that $\GC_{p_1}$ and $\GC_{p_2}$ intersect transversally at $x$ 
\end{Lem}
\begin{proof}
 We can assume $x=(0,0)$, and we use the local coordinate $(\a,\b)\in U$ 
 around $x$. 
 Each tangent vector on $T_{(0,0)}M$ is 
 given by $(\dot\a,\dot\b)\in\C\times\R\cong T_{(0,0)}(\C\times\R)$. 
 Notice that the tangent vector $(\dot\a,\dot\b)\in T_{(0,0)}M$ induces to 
 the vector field 
 $$ \(F_{1\, *}(\dot\a,\dot\b)(z), F_{2\, *}(\dot\a,\dot\b)(z) \) $$ 
 along $\GD_{(0,0)}$. 
 Here we identified $\C\times\C$ with the tangent vectors on each point of 
 $\C\times\C\subset Z$. 
 $F_{1\, *}(\dot\a,\dot\b)$ and $F_{2\, *}(\dot\a,\dot\b)$ 
 are holomorphic functions on $\bD$ 
 and their zeros coincide since 
 $$ \begin{aligned}
  F_{1\, *}(\dot\a,\dot\b)(e^{i\c}) &= i e^{i\c}\, \Xi_*(\dot\a,\dot\b)(\c), \\
  F_{2\, *}(\dot\a,\dot\b)(e^{i\c}) &= i e^{i\c}\, \overline{\Xi_*(\dot\a,\dot\b)(\c)}
  \end{aligned} $$
 by (\ref{eq:GF_i*}). 
 If $\dot\b\neq 0$, then $F_{1\, *}(0,\dot\b)(z)$ is not zero function since 
 $\GF_*$ is bijective, and $F_{1\, *}(0,\dot\b)(0)=0$ by definition. 
 This means that $(0,\dot\b)\in T_{(0,0)}M$ tangents to $\GC_{(0,0)}$ 
 since the one-parameter family of holomorphic disks fixing $(0,0)\in \GD \subset Z$ 
 is unique and this family corresponds to the vector field 
 $\(F_{1\, *}(0,\dot\b)(0), F_{2\, *}(0,\dot\b)(0) \)$
 along $\GD$. 
 
 Now consider the vector field 
 $$\(F_{1\, *}(t\dot\a,\dot\b)(0), F_{2\, *}(t\dot\a,\dot\b)(0) \)$$ 
 for $t\in [0,1]$ and non zero $\dot\a\in\C$ with sufficiently small $|\dot\a|$. 
 Then $F_{1\, *}(t\dot\a,\dot\b)$ is non zero holomorphic function on $\bD$ 
 for all $t$, and its zeros vary continuously depend on $t$. 
 Hence there is $z_1\in\bD$ near $0$ such that 
 $F_{1\, *}(\dot\a,\dot\b)(z_1)=0$, 
 and we obtain $z_1\neq 0$ since $F_{1\, *}(\dot\a,\dot\b)(0)=\dot\a\neq 0$. 
 If we put $p_2=\(F_1(0,0)(z_1),F_2(0,0)(z_1)\)\in \GD_{(0,0)}$, 
 then we find that $(\dot\a,\dot\b)\in T_{(0,0)}M$ tangents to $\GC_{p_2}$. 
 Hence $p_1=(0,0)$ and $p_2$ satisfies the statement. 
\end{proof}

\section{Construction of Einstein-Weyl spaces}
\label{Section:Construction_of_E-W_spaces}

In this Section, we construct an Einstein-Weyl structure on the parameter space 
of the family of holomorphic disks on $(Z,N)$ constructed in the previous Sections. 
The following proposition is critical. 

\begin{Prop}\label{prop:machine}
 Let $M$ be a smooth connected $3$-manifold, 
  $\varpi: \cX\to M$ be a smooth $\CP^1$-bundle. 
 Let $\rho : \cX\to \cX$  be an involution
 which commutes with $\varpi$, and is fiber-wise anti-holomorphic. 
 Suppose $\rho$ has a fixed-point set $\cX_\rho$ 
 which is an $S^1$-bundle over $M$, and which disconnects
 $\cX$ into two closed $2$-disk bundles $\cX_\pm$
 with common boundary $\cX_\rho$. 
 Let $\Dye \subset T_\C \cX$ be a distribution of complex $3$-planes 
 which satisfies the following properties: 
 \begin{itemize}
  \item $\rho_* \Dye = \overline{\Dye}$; 
  \item the restriction of $\Dye$ to $\cX_+$
      is $C^k$, $k\geq 1$ and  involutive; 
  \item $\Dye + \overline{\Dye} =T_\C\cX$ 
      on $\cX\setminus\cX_\rho$; 
  \item $\Dye \cap \ker \varpi_*$  is the $(0,1)$ tangent space of the 
      $\CP^1$ fibers of $\varpi$; 
  \item the restriction of $\Dye$ to a fiber of $\cX$ has 
      $c_1= -4$ with respect to the complex orientation; 
  \item the map $\cX\to \bP(TM) : z\mapsto 
      \varpi_*(\Dye\cap\overline\Dye)_z$ is not constant
      along each fiber of $\varpi$. 
 \end{itemize}
 Then $M$ admits a unique $C^{k-1}$ indefinite Einstein-Weyl structure 
 $([g],\nabla)$ such that 
 the null-surfaces are the projections via $\varpi$ of the integral manifolds of 
 real $2$-plane distribution $\Dye\cap T\cX_\rho$ on $\cX_\rho$; 
\end{Prop}

\begin{proof}
 Let $V^{0,1}$ be the $(0,1)$ tangent space of the fibers, then 
 $\mho=\Dye / V^{0,1}$ is a rank two vector bundle on $\cX$. 
 We can define a continuous map $\psi:\cX\to Gr_2(T_\C X)$ 
 by $z\mapsto \varpi_*(\Dye |_z)$ 
 which makes the following diagrams commute: 
 \begin{equation}
  \xymatrix{
   \cX \ar[rr]^{\psi\quad} \ar[rd] && Gr_2(T_\C X) \ar[ld] \\
   & X & }  
  \qquad 
  \xymatrix{
   \cX \ar[rr]^{\psi\quad} \ar[d]_\varrho && Gr_2(T_\C X) \ar[d]^c \\
   \cX \ar[rr]^{\psi\quad} && Gr_2(T_\C X) } 
 \end{equation}
 Using the involutiveness of $\Dye$, 
 we can prove that $\psi$ is fiber-wise holomorphic by the similar 
 argument as in \cite{bib:LM02,bib:LM05}. 

 Let 
 $ {\mathfrak P} : Gr_2(T_\C X) \longrightarrow \bP(\wedge^2 T_\C X) \cong 
   \bP(T^*_\C X) $
 be the natural isomorphism, then we obtain the fiber-wise holomorphic map 
 $\hat{\psi}={\mathfrak P}\action\psi: \cX\to \bP(T^*_\C X)$. 
 By definition, we obtain $\hat{\psi}^*\cO(-1)=\wedge^2\mho$. 
 On the other hand, since $c_1(V^{0,1})=-2$ and $c_1(\Dye)=-4$ 
 on any fiber of $\varpi$, we have $c_1(\wedge^2\mho)=c_1(\mho)=-2$. 
 Hence $\hat{\psi}$ is fiber-wise degree $2$. 
 For each fiber, there are only two possibilities for $\hat{\psi}$; 
 either a non-degenerate conic or a ramified double cover of a projective line 
 $\CP^1\subset\CP^2$. 
 
 The latter is, however, removable. 
 Indeed, any line $\CP^1\subset \CP^2$ corresponds to the planes in $\C^3$ 
 containing a fixed line. 
 Notice that, for each $z\in\cX\setminus\cX_\R$,
 $$ \varpi_*(\Dye \cap\overline{\Dye})_z =
    \varpi_*(\Dye |_z)\cap\varpi_*(\overline{\Dye}|_z)=
    \varpi_*(\Dye |_z)\cap\varpi_*(\Dye|_{\rho(z)}) $$ 
 is independent on $z$ 
 if the image of $\varpi^{-1}(x)$ under $\hat{\psi}$ is a line. 
 This contradicts to the hypothesis. 

 Now, we define a conformal structure $[g]$. 
 Let $U\subset M$ be an open set and 
 let $U\times\CP^1 \overset\sim\to \cX|_U$ be a trivialization on $U$. 
 Let $\z$ be an inhomogeneous coordinate on $\CP^1$ such that 
 $\rho(x,\z)=(x,\bar{\z})$. 
 Then we can choose a $C^k$ frame field $\{e_1,e_2,e_3\}$ on $TM|_U$ 
 so that 
 \begin{equation} \label{eq:description_of_psi}
  \hat{\psi}(x,\z)= \left[ (1+\z^2)e^1+ (1-\z^2) e^2 + 2\z e^3 \right]
 \end{equation}
 where $\{e^i\}$ is the dual frame. 
 Define an indefinite metric $g$ on $U$ so that 
 $g(e_i,e_j)$ is given by (\ref{eq:-++}). 
 Here, the frame $\{e_i\}$ is uniquely defined by (\ref{eq:description_of_psi}) 
 up to scalar multiplication, 
 and the coordinate change of $\z$ cause an $SO(1,2)$ 
 action on the frame $\{e_i\}$. 
 Hence the conformal structure $[g]$ is well-defined by $\hat{\psi}$. 
 So we obtained an indefinite conformal structure $[g]$ on $M$. 

 \vspace{2mm}
 Next we prove that a unique torsion-free connection $\nabla$ on $TM$ 
 is induced, and $([g],\nabla)$ gives an Einstein-Weyl structure on $M$. 
 We also prove that $\Dye$ agrees to the distribution $\zE$ defined in 
 Section \ref{Section:Einstein-Weyl_structures}. 

 We fix an indefinite metric $g\in[g]$, and take a local frame field 
 $\{e_1,e_2,e_3\}$ of $TM$ on a open set $U\subset M$ as above. 
 It is enough to construct $\nabla$ on $U$. 
 Notice that (\ref{eq:description_of_psi}) gives the natural identification 
 $\cX\overset\sim\to \cZ=\bP(N(T^*_\C M))$ on $U$. 
 If we define the map $\Gm_i:U\times\C\to TM$ by 
 \begin{equation} 
  \Gm_1=-e_1+e_2+\z e_3, \quad \Gm_2= \z(e_1+e_2) - e_3, 
 \end{equation}
 then we obtain $\varpi_*(\Dye |_{(x,\z)})=\Span\<\Gm_1,\Gm_2\>$  
 (cf.(\ref{eq:Gm_1&Gm_2})). 

 Let $\tilde\Gm_i$ be the vector fields on 
 $U\times\C\subset U\times\CP^1\simeq\cX|_U$ such that 
 $\tilde\Gm_i\in\Dye$ and $\tilde\Gm_i$ is written in the following form: 
 \begin{equation}
  \tilde{\Gm}_1=\Gm_1+\a \pd{\z}, \quad
  \tilde{\Gm}_2=\Gm_2+\b \pd{\z} 
 \end{equation}
 where $\a$ and $\b$ are functions on $\cX$. 
 Then $\a$ and $\b$ are uniquely defined and $C^k$. 
 Moreover, $\a$ and $\b$ are holomorphic for $\z$ since 
 $$ \left[\pd{\bar{\z}}\, ,\, \tilde{\Gm}_1 \right]
   =\frac{\del \a}{\del \bar{\z}} \pd{\z} 
   \equiv 0 \qquad \mod \Dye, $$
 and so on. 

 By the similar argument for $\z^{-1}\Gm_i$ on 
 $\{ (x,\z)\in U\times\CP^1 \, | \, \z\neq 0 \}$, 
 we find that $\z^{-1}\a\pd{\z}$ and $\z^{-1}\b\pd{\z}$ extends to 
 holomorphic vector fields on $\{\z\neq 0\}$, hence we can write 
 \begin{equation}
  \begin{aligned}
   \tilde{\Gm}_1&=\Gm_1 + (\a_0+\a_1\z+\a_2\z^2+\a_3\z^3)\pd{\z}, \\
   \tilde{\Gm}_2&=\Gm_2 + (\b_0+\b_1\z+\b_2\z^2+\b_3\z^3)\pd{\z},
  \end{aligned}
 \end{equation}
 where $\a_i$ and $\b_i$ are $C^k$ functions on $U$. 

 Recall that the compatibility condition $\nabla g=a\otimes g$ 
 holds if and only if the connection form $\w$ of $\nabla$ is written 
 \begin{equation} \label{eq:conn_form_again}
  \w= (\w^i_j)= \begin{pmatrix} 
       \p & \y^1_2 & \y^1_3 \\[1mm]
       \y^1_2 & \p & \y^2_3 \\[1mm]
 	  \y^1_3 & -\y^2_3 & \p \end{pmatrix} 
 \end{equation}
 with respect to the frame $\{e_i\}$ (c.f.(\ref{eq:conn_form})). 
 For each vector $v\in TU$, the horizontal lift $\tilde{v}$ with respect to 
 the connection defined from (\ref{eq:conn_form_again}) is 
 given by (\ref{horizontal_lift_of_pseudo_Riemm_case}). 
 If $\tilde{\Gm}_{i\,(x,\z)}$ is the horizontal lift of $\Gm_i(\z)_x$, 
 then $\y^i_j$ must be 
 \begin{equation}
  \y^2_3=\y^2_{3,0} + f e^1, \quad 
  \y^1_3=\y^1_{3,0} + f e^2, \quad 
  \y^1_2=\y^1_{2,0} - f e^3, 
 \end{equation}
 where $f$ is an unknown function on $U$ and 
 \begin{equation} \label{eq:etas_with_ambiguity}
  \begin{aligned}
   \y^2_{3,0} &= \frac{\a_0+\a_2+\b_1+\b_3}{2} e^1  
     + \frac{-\a_0-\a_2+\b_1+\b_3}{2} e^2 
	 + (-\a_3-\b_0) e^3, \\
   \y^1_{3,0} &= \frac{\a_0-\a_2+\b_1-\b_3}{2} e^1 
     + \frac{-\a_0+\a_2+\b_1-\b_3}{2} e^2 
     + (\a_3-\b_0) e^3, \\
   \y^1_{2,0} &= \frac{-\a_1+\a_3+\b_0-\b_2}{2} e^1 
     + \frac{\a_1+\a_3-\b_0-\b_2}{2} e^2. 
  \end{aligned}
 \end{equation}

 We claim that there is a unique pair $(f,\p)$ such that the connection 
 (\ref{eq:conn_form}) is torsion-free, i.e. $\w$ satisfies 
 \begin{equation} \label{eq:torsion-free_condition}
  de^i+\sum \w^i_j e^j=0. 
 \end{equation} 
 First, we fix a connection whose connection form is 
 $$ \w_0= (\w^i_{j,0})=\begin{pmatrix} 
       0 & \y^1_{2,0} & \y^1_{3,0} \\[1mm]
       \y^1_{2,0} & 0 & \y^2_{3,0} \\[1mm]
 	  \y^1_{3,0} & -\y^2_{3,0} & 0 \end{pmatrix}. $$
 Let $\l_i$ be the fiber coordinate on $T^*_\C X$ with respect to $\{e^i\}$. 
 We consider the distribution $\pi^*\Dye$ on $\cN=N(T^*_\C M)\setminus 0_M$ 
 where $\pi:\cN\to\cZ\simeq \cX$ is the projection. 
 We define 1-forms $\c, \c_{i,0}, \t_{ij,0}$ on $\cN$ (c.f.(\ref{eq:c&c_i&t_ij})) by 
 $$ \c=\sum\l_i e^i, \quad \c_{i,0}=d\l_i-\sum \l_j\w^j_{i,0}, \quad 
  \t_{ij,0}=\l_i \c_{j,0} - \l_j \c_{i,0}. $$
 If we simply write  $\t=\t_{23,0}$, then we have 
 (c.f.(\ref{eq:t_23_of_indefinite_case}))
 $$ \t= \l_2d\l_3- \l_3d\l_2 
   - \l_1 \(\l_1\y^2_{3,0} + \l_2\y^1_{3,0} - \l_3\y^1_{2,0}  \). $$
 Similar to the proof of Proposition \ref{prop:integrability_of_complex_case} or 
 \ref{prop:integrability_of_Riemannian_case}, we obtain that 
 $\pi^*\Dye=\{v\in T\cN \, | \, \c(v)=\t_{ij,0}(\c)=0\}$. 
 Hence 1-forms $\{\c,\t_{ij}\}$ are involutive. 
 
 Since $\sum\c_{i,0}\wedge e^i \equiv 0 \mod \<\c,\t_{ij}\>$, 
 we obtain $d\c \equiv \mu \mod \<\c,\t_{ij}\>$ where 
 \begin{equation} \label{eq:mu}
   \mu= (\l_2\y^1_{2,0}+\l_3\y^1_{3,0})\wedge e^1 
         + (\l_1\y^1_{2,0}-\l_3\y^2_{3,0})\wedge e^2 
		 + (\l_1\y^1_{3,0}+\l_2\y^2_{3,0})\wedge e^3 
		 +\sum \l_i de^i 
 \end{equation}
 Then we can write
 \begin{equation}
  \mu=\mu_{23}e^2\wedge e^3 + \mu_{31}e^3\wedge e^1 
     + \mu_{12}e^1\wedge e^2, 
 \end{equation}
 where $\mu_{ij}=\mu_{ij}^{\ \ l} \l_l$ are linear in $\l$. 
 Notice that  $\mu_{ij}^{\ \ l}$ are $C^{k-1}$ functions since $\c$ is $C^k$. 
 Since $d\c\equiv 0 \mod \<\c,\t_{ij}\>$, there are 1-forms 
 $\Theta_1$ and $\Theta_2$ such that 
 \begin{equation} 
    \mu= \Theta_1\wedge \t + \Theta_2\wedge \c. 
 \end{equation}
 $\Theta_1$ is, however, zero 
 since $\mu$ does not contain $d\l_i$. 
 Hence we obtain $\mu\wedge\c=0$, and this is equivalent to 
 \begin{equation}
  \begin{array}{c}
  -\mu_{23}^{\ \ 1}=\mu_{31}^{\ \ 2} = \mu_{12}^{\ \ 3}, \\[1mm]
  \mu_{12}^{\ \ 2}+\mu_{31}^{\ \ 3}=0, \quad 
  \mu_{23}^{\ \ 3}+\mu_{12}^{\ \ 1}=0, \quad
  \mu_{31}^{\ \ 1}+\mu_{23}^{\ \ 2}=0.
 \end{array}
 \end{equation}
 Thus, if we put $f=\frac{1}{2} \mu_{12}^{\ \ 3}$ and 
 $ \p=\mu_{31}^{\ \ 3}e^1+\mu_{12}^{\ \ 1}e^2+\mu_{23}^{\ \ 2}e^3$, 
 then 
 $$ \mu=-\p\wedge\c + f \(-\l_1 e^2\wedge e^3 
      + \l_2 e^3\wedge e^1+ \l_3 e^1\wedge e^2\). $$
 Here $f$ and $\p$ are $C^{k-1}$. 
 Comparing the coefficients of $\l_i$ with (\ref{eq:mu}), we obtain 
 $$ \begin{aligned}
 de^1+ \p\wedge e^1 + (\y^1_{2,0}-fe^3)\wedge e^2 
   + (\y^1_{3,0}+fe^1)\wedge e^3 &= 0, \\
 de^2 + (\y^1_{2,0}-fe^3)\wedge e^1 + \p\wedge e^2 
   + (\y^2_{3,0}+fe^1)\wedge e^3 &= 0, \\
 de^3 + (\y^1_{3,0}+fe^1)\wedge e^1 - (\y^2_{3,0}+fe^1) \wedge e^2  
   + \p\wedge e^3 &= 0. 
 \end{aligned} $$
 These are nothing but the torsion-free condition for 
 the connection defined from $f$ and $\p$ above. 

 Since $(f,\p)$ is uniquely defined, 
 we have obtained the unique torsion-free $C^{k-1}$ connection $\nabla$. 
 For this $\nabla$, the distribution $\zE$ on $\cZ\simeq \cX$ 
 agrees to $\Dye$ from the construction. 
 Hence $([g],\nabla)$ is Einstein-Weyl 
 from Proposition \ref{prop:integrability_of_p-Riemannian_case}. 
 The rest condition is deduced from the fact that 
 $\Dye\cap T\cX_\r$ corresponds to $\zD_\R$. 
\end{proof}

\begin{Rem}
 In the statement of Proposition \ref{prop:machine}, the last hypothesis
 \begin{itemize}
  \item  $\varpi_*(\Dye\cap\overline\Dye)_z$ is not constant along the fiber 
 \end{itemize}
 is not removable. 
 Actually, $\varpi_*(\Dye\cap\overline\Dye)_z$ can be constant 
 when the metric degenerates so that the light cone degenerates to a line, 
 which occurs as a limit of indefinite metric. 
\end{Rem}

\begin{Prop} \label{prop:construction_of_EW_str_on_para_space}
 Let $N$ be any embedding of $\CP^1$ into $Z=\CP^1\times\CP^1$ which 
 is $C^{2k+5}$ close to the standard one. 
 Let $\{\GD_x\}_{x\in S^2\times \R}$ be the constructed family of closed 
 holomorphic disks on $(Z,N)$. 
 Then a $C^k$ indefinite Einstein-Weyl structure $([g],\nabla)$ 
 is naturally induced on $M=S^2\times\R$. 
\end{Prop}
\begin{proof}
 We apply Proposition \ref{Prop:family_of_perturbed_holo_disks} by putting 
 $k+3$ instead of $k$ and $l=k+2$. 
 Let $M\overset\varpi\leftarrow \cX_+\overset\Gf\to Z$ be 
 the constructed double fibration (the diagram (\ref{eq:perturbed_double_fibration})), 
 then $\Gf$ is $C^{k+2}$ in this case. 
 Let $\cX_-$ be the copy of $\cX_+$ and let 
 $\cX=\cX_+\cup \overline{\cX}_-$ be the 
 $\CP^1$ bundle over $X$ which is obtained by identifying the boundaries 
 $\del \cX_+$ and $\del \cX_-$ where $\overline{\cX}_-$ is the copy of 
 $\cX_-$ with fiber-wise opposite complex structure. 
 Let $\rho : \cX\to\cX$ be the involution which interchanges 
 $\cX_+$ and $\cX_-$. 

 Let $\Gf_*:T_\C X\to T_\C Z$ be the differential of $\Gf$. 
 We define $\Dye=\Gf_*^{-1}(T^{0,1}Z)$ on $\cX_+$. 
 Then, along $\cX_\R=\del\cX_+$, $\Dye$ is spanned by $\pd{\bar{\z}}$ and the 
 distribution of real planes tangent to the fibers of 
 $\Gf: \cX_\R \to N$. 
 So we can extend $\Dye$ to whole $\cX$ so that
 $\Dye=\rho^*\overline{\Dye}$ on $\cX_\R$. 
 Let us check the hypotheses in Proposition \ref{prop:machine}. 
 \begin{itemize}
  \item $\rho_*\Dye=\overline{\Dye}$ follows from the construction. 
  \item $\Dye$ is $C^{k+1}$ on $\cX_+\setminus\cX_\R$ since $\Gf_*$ is $C^{k+1}$, 
     and $\Dye$ is involutive since $T^{0,1}Z$ is involutive. 
  \item $\Dye+\overline{\Dye}=\Gf_*^{-1}(T^{0,1}Z)+\Gf_*^{-1}(T^{1,0}Z)=
      \Gf_*^{-1}(T_\C Z)=T_\C \cX_+$ on 
      $\cX_+\setminus\cX_\R$ since $\Gf$ is surjective. 
  \item For each fiber $\varpi^{-1}(x)=\cX_+|_x$, 
      the restriction $\Gf_x:\cX_+|_x \to Z$ of $\Gf$ is a holomorphic embedding. 
	  Hence $\Dye\cap\ker\varpi_*=(\Gf_x)_*^{-1}(T^{0,1}Z)=V^{0,1}$. 
  \item $\Dye$ is $C^0$-close to the one of the standard case, 
      so $c_1(\Dye)=-4$ on each fiber of $\varpi$. 
  \item For each $x\in M$, there are $p,q\in \GD_x$ such that $\GC_p$ and $\GC_q$ 
      intersects transversally at $x$ 
	  (Lemma \ref{lem:non-degeneracy_of_time-like_vectors}). 
      If we put $z=\Gf_x^{-1}(p)=\Gf^{-1}(p)\cap\varpi^{-1}(x)$, then we obtain 
	  $$ (T_x\GC_p)\otimes\C= \varpi_*(T_{\C\, z}\Gf^{-1}(p)) 
	     = \varpi_*(\ker\Gf_*)_z = \varpi_*(\Dye\cap\overline{\Dye})_z. $$
	  Similarly $(T_x\GC_q)\otimes\C= \varpi_*(\Dye\cap\overline{\Dye})_{z'}$ 
	  for $z'=\Gf^{-1}_x(q)$. 
	  Hence $\varpi_*(\Dye\cap\overline{\Dye})$ is not constant. 
 \end{itemize}
 Thus all the hypotheses in Proposition \ref{prop:machine} are fulfilled, 
 so we obtain the unique $C^k$ indefinite Einstein-Weyl structure on $M$. 
\end{proof}

Recall that we have obtained a lift $\tilde\Gf:\cX_\R\to S(TN)$ 
of $\Gf:\cX_\R\to N$ in Section \ref{Section:The_double_fibration}. 

\begin{Prop} \label{Prop:distributions_from_Gf}
Identifying $\cX$ with $\cZ$,  
 \begin{enumerate}
  \item $\zE=\Gf_*^{-1}(T^{0,1}Z)$ on $\cX_+$ where 
    $\Gf_* : T_\C\cX_+ \to T_\C Z$, 
  \item $\zD_\R=\zE\cap T\cX_\R=\ker\{\Gf_*: T\cX_\R \to TN\}$ 
    on $\cX_\R$, 
  \item $L=\ker\{\Gf_*: T\cX_+ \to TZ\}$ on $\cX_+\setminus\cX_\R$, and 
  \item $L=\ker\{\tilde\Gf_* : T\cX_\R\to S(TN)\}$ on $\cX_\R$.
 \end{enumerate}
\end{Prop}
\begin{proof}
 {\it 1} and {\it 2} follows from Proposition \ref{prop:machine},  
 \ref{prop:construction_of_EW_str_on_para_space} and their proof. 
 We also have $\overline\zE=\Gf_*^{-1}(T^{1,0}Z)$, so 
 $L\otimes\C=\zE\cap\overline\zE=\ker{\Gf_*:T_\C\cX_+ \to T_\C Z}$. 
 Hence {\it 3} follows. 
 
 \vspace{2mm}
 Let us prove {\it 4}. Let $U\times\CP^1\overset\sim\to\cX|_U$ be a trivialization 
 on $U$ such that $\rho(x,\z)=(x,\bar{\z})$. 
 Notice that $\cX_\pm|_U=\{(x,\z)\in U\times\CP^1\, | \, \pm\Imag\z\geq 0\}$. 

 Let us denote $\z=\x+\sqrt{-1}\y$ using a real coordinate $(\x,\y)$. 
 We fix a point $(x_0,\x_0)\in\cX_\R|_U$ and let $c(s)$ be a curve defined by 
 $I_\e\to\varpi^{-1}(x) : s\mapsto (x_0,\x_0+\sqrt{-1}s)$ where 
 $I_\e=(-\e,\e)$ is a small interval. 
 Now, we define a map $\Phi:I_\e\times I_\e\to \cX: (s,t)\mapsto\Phi(s,t)$ so that 
 $\Phi(s,0)=c(s)$ and $\Phi_*(\pd{t})=l^\dagger$ where $l^\dagger$ is a 
 $\rho$-invariant real vector field such that $L=\Span\<l^\dagger\>$. 
 
 Let $\Sigma$ be the image of $\Phi$, 
 and let $\nu=\Phi(\pd{s})$ which is a tangent vector field along $\Sigma$ 
 such that $T\Sigma=\Span\<l^\dagger, \nu\>$. 
 Moreover, $\nu$ is proportional to $\pd{\y}$ on $\Sigma\cap\cX_\R$. 
 Indeed, we have $\rho\action\Phi(s,t)=\Phi(-s,t)$ by definition, 
 so $\rho_*\nu=-\nu$. 
 Hence $\nu$ is ``pure imaginary" on $\cX_\R$, i.e. 
 we can write $\nu=a\pd{\y}$ using a real-valued function $a$ on $\cX_\R$. 
 Taking $\e$ small, we can assume $a$ is positive since 
 $\nu_{(x_0,\x_0)}=c_*(\pd{s})=\pd{\y}$.
 
 Since $\{l^\dagger,\nu\}$ is involutive, there are functions $A,B$ 
 on $\Sigma$ such that $[l^\dagger,\nu]=Al^\dagger+B\nu$. 
 Let $\vp$ be a positive function on $\Sigma$ such that 
 $l^\dagger\vp=-B$, then $[l^\dagger,\vp\nu]=\vp Al^\dagger$.
 We define a positive function $\psi$ on $\Sigma\cap\cX_\R$ by 
 $\vp\nu=\psi\pd{\y}$. 
 
 Now, $\Gf:\cX_+\to Z=\CP^1\times\CP^1$ is described as 
 $\Gf(x,\z)=(F_1(x,\z),F_2(x,\z))$ on the neighborhood of $(x_0,\x_0)$ 
 using functions $F_i$ which are holomorphic on $\z$. 
 Let $p_1:Z\to\CP^1$ be the first projection, then 
 its restriction $p_1:N\to\CP^1$ is diffeomorphism. 
 Hence, identifying $N$ with $\CP^1$ by $p_1$, 
 $\Gf:\cX_\R\to N$ is described by $F_1$. 
 Since $L=\Span\<l^\dagger\>=\ker\Gf_*$ on $\cX_+\setminus\cX_\R$, 
 we have $l^\dagger F_i$=0 on $\cX_+$. Then 
 $$ l^\dagger(\vp\nu F_i)=[l^\dagger, \vp\nu] F_i + \vp\nu(l^\dagger F_i)=0, $$
 hence $l^\dagger\(\psi\frac{\del F_i}{\del\y}\) =0$ on $\Sigma\cap\cX_\R$. 
 
 Since $F_i$ are holomorphic for $\z$, we have 
 $ \frac{\del F_i}{\del \x}=-\sqrt{-1}\frac{\del F_i}{\del \y}. $
 Thus we have obtained 
 \begin{equation} \label{eq:l^dagger_is_vertical_for_tilde_Gf}
  l^\dagger\(\psi\frac{\del F_i}{\del\x}\) =0 
 \end{equation}
 on $\Sigma\cap\cX_\R$. 
 Since $\tilde{\Gf}(x,\x)=\left[\frac{\del F_1}{\del\x}(x,\x)\right]$ by definition, 
 and since $\psi$ is positive function, 
 (\ref{eq:l^dagger_is_vertical_for_tilde_Gf}) means 
 $\tilde\Gf_*(l^\dagger)=0$. 
 From {\it 2} of Proposition \ref{Prop:properties_over_S_t_general_version}, 
 the fiber of $\tilde\Gf$ is at most one-dimensional, 
 hence $L=\ker\{\tilde\Gf_*:T\cX_\R\to S(TN)\}$ on $\cX_\R$. 
\end{proof}

\begin{Prop} \label{prop:geometry_of_the_constructed_EW_space}
 The Einstein-Weyl structure $([g],\nabla)$ constructed in 
 Proposition \ref{prop:construction_of_EW_str_on_para_space}
 satisfies the following properties.  
 \begin{enumerate}
  \item For each $p\in N$, 
    $\GS_p= \left\{x\in M \, | \, p\in \del \GD_x \right\}$ 
	is connected maximal null surface 
    on $M$ and every null surface can be written in this form. 
  \item For each $p\in Z\setminus N$, 
    $\GC_p=\left\{x\in M \, | \, p\in \GD_x \right\}$ 
    is connected maximal time-like geodesic and every 
	time-like geodesic on $M$ can be written in this form. 
  \item For each $p\in N$ and non-zero $v\in T_pN$, 
    $\GC_{p,v}= \left\{x\in M \, | \, p\in \del \GD_x, v\parallel \GD_x \right\}$ is 
	connected maximal null geodesic on $M$ and every 
	null geodesic on $M$ can be written in this form. 
 \end{enumerate}
\end{Prop}
\begin{proof}
 From Proposition \ref{Prop:distributions_from_Gf} and 
 the properties of $\zD_\R$ and $L$, we obtain 
 \begin{itemize} 
  \item $\GS_p =\varpi\action\Gf^{-1}(p) $ is a null surface for each $p\in N$, 
  \item $\GC_p =\varpi\action\Gf^{-1}(p)$ is a time-like geodesic 
      for each $p\in Z\setminus N$, 
  \item $\GC_{p,v}=\varpi\action\tilde\Gf^{-1}([v])$ is a null geodesic 
      for each $p\in N$ and non zero $v\in T_pN$. 
 \end{itemize}
 Moreover from Proposition \ref{Prop:properties_over_S_t_general_version}, 
  \begin{itemize} 
  \item $\GS_p\simeq S^1\times\R$ for each $p\in N$, 
  \item $\GC_p\simeq \R$ for each $p\in Z\setminus N$, 
  \item $\GC_{p,v}\simeq \R$  for each $p\in N$ and non zero $v\in T_pN$, 
 \end{itemize}
 and they are all closed in $M$. Hence the statement follows. 
\end{proof}

Recall the compactification of the double fibration given by 
(\ref{eq:compactified_double_fibration_perturbed_case}). 
Let $\hat{\GC}_p$ and $\hat\GC_{p,v}$ be the compactification of 
$\GC_p$ and $\GC_{p,v}$ in $\hat{\cX}_+$ respectively.  

\begin{Prop} \label{prop:foliation_on_N}
 \begin{enumerate}
  \item For each $p\in Z\setminus N$, 
    $\hat{\cX}_\R|_{\hat{\GC}_p}$ is homeomorphic to $S^2$ and the restriction 
	$\hat{\Gf}:\hat{\cX}_\R|_{\hat{\GC}_p}\to N$ is a homeomorphism.  
    In particular, $\{\del\GD_x\}_{x\in\GC_p}$ gives a foliation on 
	$N\setminus \{\text{\rm 2 points}\}$.  
  \item For each $p\in N$ and non zero $v\in T_pN$, 
    $\hat{\cX}_\R|_{\hat{\GC}_{p,v}}$ is homeomorphic to $S^2$ and the restriction 
	$\hat{\Gf}: \hat{\cX}_\R|_{\hat{\GC}_{p,v}} \to N$ is surjective. 
	Moreover, this is one-to-one distant from the curve $\hat{\Gf}^{-1}(p)$, hence 
    $\{(\del\GD_x\setminus\{p\})\}_{x\in\GC_{p,v}}$ gives a foliation on 
	$N\setminus \{p\}$.  
 \end{enumerate}
\end{Prop}
\begin{proof}
 Let $p\in Z\setminus N$, then $\cX_\R|_{\GC_p}$ is an $S^1$-bundle over 
 $\GC_p\simeq\R$. Since $\hat{\cX}_\R|_{\hat{\GC}_p}$ is the compactification 
 of $\cX_\R|_{\GC_p}$ with extra two points, it is isomorphic to $S^2$. 
 Since $\Gf$ is $C^0$-close to the one of the standard case, 
 $\hat{\Gf}: \hat{\cX}_\R|_{\hat{\GC}_p} \to N$ is degree one map. 
 
 Let $\Gf_*:T(\cX_\R|_{\GC_p})\to TZ_\R$ be the differential. 
 We claim that $\ker \Gf_*=0$ everywhere. 
 Indeed, if there exist non zero $w\in T_z(\cX_\R|_{\GC_p})$ 
 such that $\Gf_*(w)=0$, then $w\in\Dye_z$ and $\varpi_*(w)\neq 0$. 
 Then $\varpi_*(w)$ must be null with respect to the constructed conformal structure. 
 On the other hand $\varpi_*(w)$ tangents to $\GC_p$, 
 so this is time-like. This is a contradiction. 

 Hence $\hat{\Gf}: \hat{\cX}_\R|_{\hat{\GC}_p} \to N$ is 
 locally homeomorphic degree one map, i.e. homeomorphism. 

 \vspace{2mm}
 Next, let $p\in N$. By the similar argument, 
 $\hat\cX_\R|_{\hat{\GC}_{p,v}}\simeq S^2$ and 
 $\hat{\Gf}:\hat{\cX}_\R|_{\hat{\GC}_{p,v}} \to N$ is degree one, hence surjective. 

 We claim that $\ker\{\Gf_*: T(\cX_\R|_{\GC_{p,v}}) \to TN\}=0$ on 
 $z\in \( \cX_\R|_{\GC_{p,v}} \setminus \Gf^{-1}(p) \)$. 
 Indeed, if there exists non zero $w\in T_z(\cX_\R|_{\GC_{p,v}})$ 
 such that $\Gf_*(w)=0$, then $\varpi_*(w)$ is non zero and null. 
 Notice that $\varpi_*(w)$ tangents to the null surface $\GS_{\Gf(z)}$. 
 
 On the other hand, $\varpi_*(w)$ tangents to $\GC_{p,v}\subset\GS_p$. 
 Since $\Gf(z)\neq p$, $\GS_{\Gf(z)}$ and $\GS_p$ are different null surfaces, 
 hence $T_{\varpi(z)}\GS_{\Gf(z)}$ and $T_{\varpi(z)}\GS_p$ are 
 different null planes at $\varpi(z)$. 
 Then $\varpi_*(w)\in T_{\varpi(z)}\GS_{\Gf(z)}\cap T_{\varpi(z)}\GS_p$ 
 must be space-like vector, this is a contradiction. 
 Hence the statement follows. 
\end{proof}

\begin{Prop} \label{prop:space-like_zoll}
 Let $([g],\nabla)$ be the Einstein-Weyl structure constructed in 
 Proposition \ref{prop:construction_of_EW_str_on_para_space}. 
 Then, for each distinguished $p,q\in N$, 
 $\GC_{p,q}= \left\{x\in M \, | \, p, q \in \del \GD_x \right\}$ is 
 connected closed space-like geodesic on $M$ and every 
 space-like geodesic on $M$ can be written in this form. 
 In particular, this Einstein-Weyl structure is space-like Zoll. 
\end{Prop}

\begin{proof}
 Since $\GC_{p,q}$ is the intersection of the null surfaces $\GS_p$ and $\GS_q$, 
 this is either empty or a space-like geodesic. 
 We claim that $\GC_{p,q}$ is not empty and is homeomorphic to $S^1$. 
 For each non zero $v\in T_pN$, there is a unique $x\in\GC_{p,v}$ such that 
 $q\in\del\GD_x$ since $\{(\del \GD_x\setminus\{p\})\}_{x\in\GC_{p,v}}$ foliates 
 $N\setminus\{p\}$ by {\it 2} of Proposition \ref{prop:foliation_on_N}. 
 Then $x\in\GC_{p,q}$, so $\GC_{p,q}$ is not empty. 
 Moreover there is a one-to-one continuous map 
 $S(T_pN)\to\GC_{p,q}$, so $\GC_{p,q}\simeq S^1$. 
\end{proof}

The main theorem (Theorem \ref{thm:main_theorem}) follows from 
Proposition \ref{prop:construction_of_EW_str_on_para_space}, 
\ref{prop:geometry_of_the_constructed_EW_space} and 
\ref{prop:space-like_zoll}.


\vspace{13mm}
\noindent
\small
\begin{tabular}{l}
Department of Mathematics \\
Graduate School of Science and Engineering \\
Tokyo Institute of Technology \\
2-12-1, O-okayama, Meguro, 152-8551, JAPAN \\
{\tt {nakata@math.titech.ac.jp}}
\end{tabular}

\end{document}